%% file: paper_Yangian.tex
\documentclass[a4paper,10pt]{article}

\input{preamble/usepackage}
\input{preamble/newcommand}
\input{preamble/texmacs}
\usepackage{mathtools}
\usepackage{cite}

\hypersetup{
    colorlinks=true,
    linkcolor=blue,
    filecolor=magenta,    
    urlcolor=cyan,
    citecolor=cyan,}
\newcommand{\g}{\mathfrak g}

\title{Stokes Phenomenon and Yangians}
\author{Qian Tang and Xiaomeng Xu}
\date{\today}

\begin{document}

\maketitle
\abstract{
In this paper, we first establish a connection between Yangians and the unique formal solution of the quantum hypergeometric differential equations at irregular singularities.
We then realize the Stokes matrices of the hypergeometric equations as infinite matrix products of representations of Yangains, with the help of the theory of difference systems. Along the way, we also investigate the algebroid structure associated with the Stokes matrices. 

}



    

\section{Introduction}


Let us take the complex Lie algebra $\mathfrak{gl}_\nu$ and its universal enveloping algebra $U(\mathfrak{gl}_\nu)$ generated by $\{e_{ij}\}_{1\leqslant i,j\leqslant n}$ subjected to the relation $[e_{ij},e_{kl}]=\delta_{kj}e_{il}-\delta_{li} e_{kj}$. 
Let us 
consider the quantum confluent hypergeometric system
\begin{equation}\label{introqeq}
\frac{\mathd}{\mathd z} F(z) =
\left( u - \hbar \mathbf{E} z^{-1}\right)\cdot F(z),
\end{equation}
where
the $\nu\times \nu$ matrix $\mathbf{E}\assign(e_{ij})_{\nu\times \nu}$ 
with entries valued in $U(\mathfrak{gl}_\nu)$,
$u=\diag(u_1,\ldots,u_\nu)$ 
is seen as an $\nu\times \nu$ matrix with scalar entries in $U(\mathfrak{gl}_\nu)$,
and $\hbar$ is a complex parameter.
For most of the time, we will take a finite-dimensional representation $V$ of $\mathfrak{gl}_\nu$, and consider the system \eqref{introqeq} for 
an $\nu\times \nu$ matrix function $F^V(z)$ with entries in ${\rm End}(V)$, i.e., a (block matrix) solution $F^V(z)\in {\rm End}(V)\otimes {\rm End}(\mathbb{C}^{\nu})$. 
Here, the action of the coefficient matrix of \eqref{introqeq} on $F^V(z)$ is given by matrix multiplication and the representation $V$ of $\mathfrak{gl}_\nu$. 

Confluent hypergeometric systems have played important roles in many subjects, \emph{e.g.} the study of Frobenius manifolds \cite{Dubrovin} and in particular the quantum cohomology of Fano manifolds \cite{GGI}, linearization in Poisson geometry \cite{Boalch}, quantum Weyl group actions on Poisson groups \cite{BoalchG}, stability conditions \cite{Br, BTL}, long time asymptotics of some isomonodromy deformation equations \cite{TX2} \emph{etc}. The system \eqref{introqeq} is a natural quantum analog of the confluent hypergeometric system. Its Stokes phenomenon and the WKB approximation have been used to give a transcendental realization of the quantum group $U_q(\mathfrak{gl}_\nu)$ and its crystal structure in \cite{Xu2, Xu}. Therefore, the system is worthy of further study. 

The purpose of this paper is then twofold; one is to establish a new relation between the system \eqref{introqeq} and the Yangian; the other is to deepen the relation between the Stokes phenomenon of \eqref{introqeq} and quantum groups, using Darboux's method and a path algebroid arising from the resurgence theory.


\subsection{Formal Solutions and Yangian \texorpdfstring{$Y(\mathfrak{gl}_{\nu-1})$}{Y(gl(ν-1))}}

For a positive integer $m$, the Yangian $Y_{\hbar}(\mathfrak{gl}_{m})$ was introduced in the formulation of quantum inverse scattering methold by Faddeev’s school (see e.g.
\cite{takhtadzhan1979quantum,
kulish1982solutions,
faddeev1982integrable})
and by Drinfeld \cite{Drinfeld1990} in the study of the Yang-Baxter equation.
It is a Hopf algebra that can be regarded as a deformation of the enveloping algebra $U(\mathfrak{gl}_m[x])$, where $\mathfrak{gl}_m[x]$ is the Lie algebra of $\mathfrak{gl}_m$-valued polynomials
\cite{molev2002yangians}.

\begin{defi}\cite{Drinfeld1990}
The Yangian for $\mathfrak{gl}_m$
is a unital associative algebra with countably
many generators $\{t_{ij}^{(p)}: p\in\mathbb{N}^+, i,j=1,\ldots,m\}$, and the generating relations for all $i_0, i_1, j_0, j_1\in \{1,...,m\}$,
\begin{equation}\label{Yangdef}
\frac{1}{\hbar}
[T(\lambda_1)_{i_1 i_0},
 T(\lambda_2)_{j_1 j_0}] =
\frac{T(\lambda_2)_{j_1 i_0}
	  T(\lambda_1)_{i_1 j_0}-
	  T(\lambda_1)_{j_1 i_0}
	  T(\lambda_2)_{i_1 j_0}}
{\lambda_1-\lambda_2}.
\end{equation}
Here $T(\lambda)_{ij}$ is the generating series $T(\lambda)_{ij} 
    = \delta_{ij} +
    \sum_{p=1}^\infty t_{ij}^{(p)} \lambda^{-p}.$
\end{defi}

\begin{thm}\label{MainThm}
For $\hbar\notin \mathbb{Q}$, the system \eqref{introqeq} has a unique formal solution
\begin{align}\label{SolA}
\hat{F}(z)= \left(I+\sum_{p=1}^\infty {H}_p z^{-p}\right)\cdot 
\e^{u z}
z^{-\hbar\delta\mathbf{E}},\quad
\delta\mathbf{E}=\diag(e_{11},\ldots,e_{\nu\nu}),
\end{align}
where $I$ is the identity matrix, and the coefficients $H_p\in\operatorname{Mat}_{n}
(U(\mathfrak{gl}_\nu{}))$ are recursively determined by the relations
\begin{align}\label{Form:QHp}
    ({H}_{1})_{i k} = 
    -\frac{\hbar e_{ik}}{u_k-u_i}
    ,\quad
    ({H}_{p+1})_{i k} = 
    \sum_{j\neq k}
    T_k(p)_{ij}
    ({H}_{p})_{j k}
    ,\quad
    ({H}_{p})_{k k} =
    \sum_{j\neq k}
    \frac{\hbar e_{k j}}{p} 
    ({H}_{p})_{j k}.
\end{align}
Here for any $\lambda\in \mathbb{C}\setminus\{0\}$, 
$T_k(\lambda)$ is an $(\nu-1)\times(\nu-1)$-matrix with entries
\begin{align}\label{Yn-1}
    T_k(\lambda)_{ij}:=
    \frac{1}{u_k-u_i}
    \left(
    (\lambda+\hbar(e_{kk}+1))\delta_{ij}
    - \hbar e_{ij}
    - \hbar^2\frac{e_{ik}e_{kj}}{\lambda}
    \right)\in U(\mathfrak{gl}_\nu),\quad
    i,j\in \{ 1,\ldots,\nu\}\setminus \{k\}.
\end{align}
Furthermore, $T_k(\lambda)$ satisfies
the Yangian relations \eqref{Yangdef} of $Y_{\hbar}(\mathfrak{gl}_{\nu-1})$
, i.e., 
for indices $i_0, i_j, j_0, j_1$,
\begin{align}\label{Intro:RTTYangian}
\frac{1}{\hbar}
[T_k(\lambda_1)_{i_1 i_0},
 T_k(\lambda_2)_{j_1 j_0}] & =
\frac{T_k(\lambda_2)_{j_1 i_0}
	  T_k(\lambda_1)_{i_1 j_0}-
	  T_k(\lambda_1)_{j_1 i_0}
	  T_k(\lambda_2)_{i_1 j_0}}
{\lambda_1-\lambda_2}\in U(\mathfrak{gl}_\nu).
\end{align}
\end{thm}

\begin{rmk}
For convenience, we always require that the $\nu-1$ indices of
the $(\nu-1)\times(\nu-1)$ matrix $T_k(\lambda)$ 
are $\{1,\ldots,k-1,k+1,\ldots,\nu\}$.
It should be emphasized that for different $\hbar \in \mathbb{C} \setminus \{0\}$, the Yangians $Y_\hbar(\mathfrak{gl}_{m})$ are canonically isomorphic.
\end{rmk}

Our first theorem states that the Yangian $Y_{\hbar}(\mathfrak{gl}_{\nu-1})$ naturally arises from the formal power series solutions of the system \eqref{introqeq}.
Equivalently, for each $k=1,\ldots,\nu$, we get an algebra homomorphism
\begin{equation}\label{Oconstruct1}
O(u)_k: Y_{\hbar}(\mathfrak{gl}_{\nu-1}) \rightarrow U(\mathfrak{gl}_\nu);\quad 
T(\lambda)_{ij}\mapsto T_k(\lambda)_{ij}, \ \ i,j=1,...,k-1,k+1,...,\nu.
\end{equation}
\begin{rmk}
Here the algebra homomorphism $O(u)_k$ is a $(u_1,...,u_\nu)$-family of deformation of the homomorphism in the Olshanski centralizer construction
\cite{olshanski1991representations,olshanski1987extension}.
\end{rmk}

\subsection{Stokes matrices as Infinite Product of Representations of Yangian}

For any fixed $u$ with distinct eigenvalues and nontrivial representation $V$, the series $I+\sum_{p=1}^\infty {H}_p z^{-p}$, in the corresponding formal solution $\hat{F}(z)$, is a divergent ${\rm End}(V)\otimes {\rm End}(\mathbb{C}^{\nu})$ valued formal power series. Thus $\hat{F}(z)$ is only a formal solution. The resummation (see, e.g., \cite{Balser}) states that
there exist certain sectorial regions around $z=\infty$, such that on each of these sectors there is a unique ${\rm End}(V)\otimes {\rm End}(\mathbb{C}^{\nu})$ valued holomorphic solution with the prescribed asymptotics $\hat{F}(z)$, see Proposition \ref{Pro:CAFS} for more details. These solutions are in general different
(that reflects the Stokes phenomenon), and the transition between them can be measured by the Stokes
matrices $S_{[\tau]}\in {\rm End}(V)\otimes {\rm End}(\mathbb{C}^{\nu})$ (see Definition \ref{monodata}), associated to the anti-Stokes direction $\tau=-\arg(u_j-u_i)$ for some $i\ne j$. For $u_1,...,u_\nu$ in generic position, the only none zero subdiagonal entry of $S_{[\tau]}$ is the $(i,j)$-entry.
Our second result realizes the Stokes matrices
as one-sided infinite matrix products of the $(\nu-1)\times (\nu-1)$ matrices $T_k(m)\in {\rm End}(V)\otimes {\rm End}(\mathbb{C}^{\nu-1})$. 

\begin{pro}\label{infprodS}
Suppose that none of $u_k$, for all $k\ne i,j$, lies on the segmant determined by $u_i$ and $u_j$, then the $(i,j)$-entry of the Stokes matrix $S_{[\tau]}$ associated to the anti-Stokes direction $\tau=-\arg(u_j-u_i)$, as elements in ${\rm End}(V)$, is given by 
(see Section \ref{Sect:DiffSys} for the meaning of the infinite matrix product)
\begin{align*}
   (u_j-u_i)^{\hbar e_{ii}}
    \frac{(S_{[\tau]})_{ij}}{2\pi\mathi}
    (u_j-u_i)^{-\hbar e_{jj}} & =
    \lim_{p\to\infty}
    \frac{(u_j-u_i)^{p+1}}{p!}
    p^{\hbar(e_{ii}-e_{jj}-1)}
    \left(
    \overleftarrow{\prod_{m=1}^{p}}
	T_j(m)
    \right)_{i\hat{j}}
    \left(
    -\frac{1}{u_jI-u_{\hat{j}\hat{j}}}
    \hbar \mathbf{E}_{\hat{j}j}
    \right),\\
    (u_j-u_i)^{\hbar e_{ii}}
    \frac{(S_{[\tau]})_{ij}}{2\pi\mathi}
    (u_j-u_i)^{-\hbar e_{jj}}
    & =
    \lim_{p\to\infty}
    \frac{(u_j-u_i)^{p}}{p!}
    \left(-\hbar \mathbf{E}_{i\hat{i}} \right)
    \left(
    \overleftarrow{\prod_{m=-p}^{-1}}
	T_i(m)
	\right)_{\hat{i} j}
  	p^{\hbar(e_{ii}-e_{jj}+1)},
\end{align*}
where $\overleftarrow{\prod_{m=1}^{p}}T_j(m) := T_j(p)\cdots T_j(2)T_j(1)$.
Here for a matrix $A$, 
we denote $A_{\hat{i}j}$ as the $j$-th column of $A$ without the $i$-th row; 
$A_{i\hat{j}}$ as the $i$-th row of matrix $A$ without the $j$-th column.
\end{pro}

Since, for each $k=1,\ldots,\nu$, the $U(\mathfrak{gl}_\nu)$-module $V$ can be seen as a representation of $Y_{\hbar}(\mathfrak{gl}_{\nu-1})$ via the algebra homomorphism $O_k: Y_{\hbar}(\mathfrak{gl}_{\nu-1}) \rightarrow U(\mathfrak{gl}_\nu)$ in \eqref{Oconstruct1}, the above proposition states that the Stokes matrices can be obtained as infinite product of the representations of Yangian.

Following \cite{Xu}, the Stokes matrices of the system \eqref{introqeq} satisfy 
the Faddeev-Reshetikhin-Takhtajan's RLL relations of quantum groups
\cite{Faddeev1988}. 
As a consequence, we deduce that the one-sided infinite matrix products of 
the $n$ representations $\{{T}_k(\lambda)\}_{k=1,\ldots,\nu}$ of $Y(\mathfrak{gl}_{\nu-1})$ on $V$ satisfy some commutative relations. That is 

\begin{thm}\label{oneside}
Let us introduce the matrix 
$\mathscr{T}$ with entries
\begin{equation}\label{infmatpro}
    \mathscr{T}_{ij} \ \assign
    \lim_{p\to\infty}
    \frac{(u_j-u_i)^{p+1}}{p!}
    p^{\hbar(e_{ii}-e_{jj}-1)}
    \left( \overleftarrow{\prod_{m=1}^{p}}
	T_j(m)
    \right)_{i\hat{j}}
    \left(
    -\frac{1}{u_jI-u_{\hat{j}\hat{j}}}
    \hbar \mathbf{E}_{\hat{j}j}
    \right),
\end{equation}
where $T_k(m)$ is defined as in \eqref{Yn-1}. Then
\begin{enumerate}
    \item For distinct indices $s_1,t_1,s_2,t_2$ and
    disjoint segments
    $[u_{s_1},u_{t_1}]\cap [u_{s_2},u_{t_2}] = \varnothing$, 
    we have
    \begin{align}\label{Case:DisjSeg}
    \mathscr{T}_{s_1 t_1} 
    \mathscr{T}_{s_2 t_2} -
    \mathscr{T}_{s_2 t_2} 
    \mathscr{T}_{s_1 t_1} = 0;
    \end{align}
    
    \item For distinct indices $s_1,t_1,s_2,t_2$ and disjoint segments
    $[u_{s_1},u_{t_1}]\cap [u_{s_2},u_{t_2}] \neq \varnothing$,
    $\arg(u_{t_2}-u_{s_2}) -
    \arg(u_{t_1}-u_{s_1}) \in (0,\pi)$, we have
    \begin{align}\label{Case:CrossSeg}
    \frac{
    \mathscr{T}_{s_1 t_1} 
    \mathscr{T}_{s_2 t_2} -
    \mathscr{T}_{s_2 t_2} 
    \mathscr{T}_{s_1 t_1}
    }{q-q^{-1}} =
    \frac{
    (u_{t_2}-u_{s_1})^{\hbar( e_{t_2t_2}-e_{s_1s_1}+1 )}
    (u_{t_1}-u_{s_2})^{\hbar( e_{t_1t_1}-e_{s_2s_2}+1 )}
    }{
    (u_{t_1}-u_{s_1})^{\hbar( e_{t_1t_1}-e_{s_1s_1}+1 )}
    (u_{t_2}-u_{s_2})^{\hbar( e_{t_2t_2}-e_{s_2s_2}+1 )}
    }
    \mathscr{T}_{s_1 t_2}
    \mathscr{T}_{s_2 t_1};
    \end{align}

    \item For distinct indices $s,m,t$, we have
    \begin{subequations}
    \begin{align}\label{Case:smt}
    \frac{
    \left(\frac
    {u_s-u_m}{u_t-u_s}
    \right)^\hbar
    \mathscr{T}_{sm} \mathscr{T}_{mt} -
    \left(\frac
    {u_t-u_m}{u_t-u_s}
    \right)^\hbar
    \mathscr{T}_{mt} \mathscr{T}_{sm}
    }{(q-q^{-1})/2\pi\mathi} & =
    -
    \left(\frac
    {u_m-u_s}{u_t-u_s}
    \right)^{\hbar(e_{ss}-e_{mm})}
    \left(\frac
    {u_t-u_m}{u_t-u_s}
    \right)^{\hbar(e_{mm}-e_{tt})}
    \mathscr{T}_{st},\\
    \mathscr{T}_{ms} \mathscr{T}_{mt} -
    \left(\frac
    {u_m-u_s}{u_t-u_m}
    \right)^\hbar
    \mathscr{T}_{mt} \mathscr{T}_{ms}
    & = 0,\\
    \mathscr{T}_{sm} \mathscr{T}_{tm} -
    \left(\frac
    {u_m-u_s}{u_t-u_m}
    \right)^\hbar
    \mathscr{T}_{sm} \mathscr{T}_{tm}
    & = 0,\\
    \label{Case:stts}
    \frac{
    \mathscr{T}_{st} \mathscr{T}_{ts} -
    \mathscr{T}_{ts} \mathscr{T}_{st}
    }{(q-q^{-1})/2\pi\mathi} & = 
    \frac{
    q^{e_{ss}-e_{tt}} -
    q^{e_{tt}-e_{ss}}
    }{2\pi\mathi}.
    \end{align}
    \end{subequations}  
\end{enumerate}
We require that for arbitrary indices 
$s_1,t_1,s_2,t_2$ (not necessarily distinct),
the arguments satisfy
$\arg(u_{t_2}-u_{s_2}) -
 \arg(u_{t_1}-u_{s_1}) \in (-\pi,\pi)$.
\end{thm}

The commutative relation of bilateral infinite matrix product, associated to one representation, of Yangians was studied in the development of the $R$-matrix formulation of the quantum inverse scattering method, see e.g., Faddeev-Reshetikhin \cite{FR} for the semiclassical limit case. In the formulation, for any fixed $k=1,...,\nu$, the Yangian relation \eqref{Intro:RTTYangian} can be used to compute the commutative relation of the bilateral infinite matrix product 
\[\lim_{p\to\infty}
N_1(p)\left( \overleftarrow{\prod_{-p}^{p}}
	T_k(m)
    \right) N_2(p)\] 
(with $N_1(p), N_2(p)$ certain normalizers). Our Theorem \ref{oneside} can be seen as a generalization of this formulation to the case of multiple (related) Yangian representations: $n$ representations $T_1,...,T_n$ are given, and for different indices $k$ and $j$, the commutator between $T_k$ and $T_j$ satisfies a quadratic relation similar to the Yangian relation \eqref{Intro:RTTYangian}. One can imagine that, these Yangian-like quadractic relations
permit calculation of the commutators between the elements of the one-sided infinite matrix products showing to have forms \eqref{Case:DisjSeg}-\eqref{Case:stts}. Besides, given Proposition \ref{infprodS}, one obtains a new interpretation of the quantum group relations of the Stokes matrices from this viewpoint. We hope to further develop this generalized formulation in future.

Let us mention that the system \eqref{introqeq} can be seen as a Knizhnik–Zamolodchikov (KZ) type equation with an irregular singularity. 
The KZ equations with irregular singularities have been introduced from various perspectives. For example, it was introduced in \cite{reshetikhin1992knizhnik}, and was given a
representation-theoretic interpretation in \cite{felder2023singular}. From the perspective of isomonodromy deformation, the works \cite{Rem, Rem2} construct an equation via quantisation of the
Hamiltonian systems of \cite{Boalch5}, and that equation reduces to KZ equation in the logarithmic
case-so it is also an irregular version of KZ equation. 
These resulting differential equations can have arbitrary order pole and therefore have Stokes phenomenon. Theorem \ref{MainThm} states that the algebraic structure hidden behind the recursive relations of the formal solutions of KZ equation with a second order pole is the Yangian. It is then interesting to discover the algebraic structures hidden in the formal solutions of KZ equation with higher order poles.

\subsection{Outlook: Classical Lie Algebra Cases and Twisted Yangians}
Let us first discuss the generalization of the above results to other types. Let $\g_\nu$ denote the rank $\nu$ simple complex Lie algebra of type 
B, C, or D, i.e.,
\[\g_n=\mathfrak{o}_{2\nu+1}, \ \ \mathfrak{sp}_{2\nu}, \ \ 
 \text{or} \ \ \mathfrak{o}_{2\nu}.\]
For $-\nu\le i,j\le \nu$, let us introduce the generators of $U(\mathfrak{g}_\nu)$
\[K_{ij}=e_{ij}-\theta_{ij} e_{ji}, \]
where
\begin{equation}
\theta_{ij}=\left\{\begin{array}{lr} 1,  & \text{ in the orthogonal case} \\ {\rm sgn}(i)\cdot {\rm sgn}(j), & \text{ in the symplectic case}.
             \end{array}
\right.
\end{equation}

Given any finite-dimensional irreducible representation $L(\lambda)_{\g_\nu}$ of $\g_n$ with a highest weight $\lambda$, let us consider the quantum confluent hypergeometric type equation
\begin{eqnarray}\label{introqeqBCD}
\frac{dF}{dz}=h\Big( u+\frac{K}{z}\Big)\cdot F,
\end{eqnarray}
for $F(z)\in {\rm End}(L(\lambda)_{\g_n})\otimes {\rm End}(\mathbb{C}^{m})$ with $m=2\nu$ or $2\nu+1$. Here the $m\times m$ matrix $K=(K_{ij})$ has entries valued in $U(\g_\nu)$, and the $m\times m$ matrix $u$ is diagonal with entries 
\begin{equation}
u=\left\{\begin{array}{lr} {\rm diag}(u_\nu, \ldots, u_1, 0, u_1, \ldots, u_\nu),  & \mathfrak{s}\mathfrak{o}_{2 n + 1} \\ {\rm diag}(u_\nu, \ldots,u_1, u_1, \ldots, u_\nu),  & \mathfrak{s}\mathfrak{o}_{2 n} \text{ or } \mathfrak{s}\mathfrak{p}_{2 n}.
             \end{array}
\right.
\end{equation}
We make the assumption that $u$ is regular, i.e., $u_1,...,u_\nu$ and $0$ are distinct. Similar to Theorem \ref{MainThm}, the equation \eqref{introqeqBCD} has a recursively defined formal power series solution. We expect that the algebraic structure hidden behind the recursive relation is the corresponding twisted Yangian. We remark that for the symplectic Lie algebra case, the Stokes matrices of the equation \eqref{introqeqBCD} (with $u$ having distinct eigenvalues) are proved in \cite{Xu1} to be the $K$-matrix for the quantum symmetric pair of type C. Therefore, it is interesting to study further the relation between the analytic theory of \eqref{introqeqBCD} and the theory of various quantum algebras.

\vspace{2mm}

The organization of the paper is as follows. 
In Section \ref{ClassicalCase}, 
we introduce the basic concepts, including the Stokes matrices, in the study of confluent hypergeometric systems.
In Section \ref{Sect:DiffSys}, 
we present the classical theory of difference equations,
where the Stokes matrices are realized as canonical analytical solutions of some difference equations and are thus interpreted as infinite product of matrices. 
In Section \ref{Sect:PathAlg}, 
we discuss the algebroid structures underlying the Stokes matrices. 
The results of this section can help us obtain the commutation relations between the entries of the Stokes matrices given by infinite product of matrices.
In Section \ref{Sect:Q}, 
we establish the connection
between the Stokes matrices of system \eqref{introqeq} and Yangians. The connections between various lemmas and proposition used, and theorems proved in the paper can be summarized in the following diagram. 

\begin{center}
\begin{tikzcd}
                                         & 
\text{Lem. \ref{Lem:AnaLim}} \arrow[d] & 
\text{Prop. \ref{Prop:MD}} \arrow[r] \arrow[ld] & 
\text{{\color{blue}Lem.} \ref{Lem:BdFdnice}} \arrow[rd] \arrow[d] & 
\text{Lem. \ref{LdFdnice}} \arrow[l]     \\
\text{Prop. \ref{FormSolirr}} \arrow[rd]            & 
\text{{\color{blue}Cor.} \ref{Cor:LimCoef}} \arrow[d]             & 
\text{{\color{blue}Cor.} \ref{Weight:Seg}} \arrow[l]              & 
\text{{\color{blue}Prop.} \ref{Pro:gammajk}} \arrow[l] \arrow[d]  & 
\text{{\color{blue}Cor.} \ref{ProveWeight}} \arrow[ld]            \\
\text{Prop. \ref{FormSolreg}} \arrow[r]     & 
\text{{\color{blue}Thm.} \ref{Thm:Difeq}} \arrow[d]       & 
\text{Prop. \ref{AppPro:Sto[t]}} \arrow[u]  & 
\text{{\color{blue}Thm.} \ref{Thm:SWeight}} \arrow[d]     & 
\text{Thm. \ref{qStokes}} \arrow[ld]        \\
\text{Thm. \ref{MainThm}} \arrow[r]             & 
\text{{\color{blue}Prop.} \ref{infprodS}} \arrow[rd]           & 
\text{Prop. \ref{DiffLim}} \arrow[d] \arrow[lu] & 
\text{{\color{blue}Prop.} \ref{Pro:PathtoUq}} \arrow[ld]      & 
\text{Lem. \ref{Lem:Leb}} \arrow[l]             \\
                            &
                            & 
\text{{\color{blue}Thm.} \ref{oneside}}   &                                                                               &                                         
\end{tikzcd}
\end{center}

\section{Monodromy Data of the Confluent Hypergeometric System}
\label{ClassicalCase}

The \textbf{$n$-th confluent hypergeometric system} 
is a $n\times n$ system of the following form
\begin{eqnarray}\label{Conflu}
    \frac{\mathd}{\mathd z} {F}(z) =
    \left(u+A z^{-1}\right)\cdot {F}(z),
\end{eqnarray}
where the solution ${F}(z)$ is an $n\times n$ matrix-valued analytic function.
System \eqref{Conflu} has only two singularities, which are also the only singularities of its solution. The irregular singularity is $z=\infty$, and the regular singularity is $z=0$.

In this section,
we always impose the following conditions and notations:
\begin{itemize}
    \item $u=\diag(u_1 I_{n_1},\ldots, u_\nu I_{n_\nu})$,
    where $u_1,\ldots,u_\nu$ are distinct, 
    with multiplicity $n_1,\ldots,\nu_\nu$;
    \item the residue matrix is divided into
    $(n_1,\ldots,\nu_\nu)$-blocks,
    denoted by $A=(A_{i j})_{\nu\times\nu}$
    the blocked matrix.
\end{itemize}

Note that the system \eqref{introqeq} associated to a representation $V$ becomes a special case of the System \eqref{Conflu} with rank $n=m \nu$. Actually, assume $m={\rm dim}(V)$, then the system \eqref{introqeq} associated to $V$ is an $m\nu\times m\nu$ system, where $u=\diag(u_1 I_{m},\ldots, u_\nu I_{m})$, with distinct $u_1,\ldots,u_\nu$ and multiplicity $m,\ldots,m$, and
 the residue $A=-\hbar\mathbf{E}^V=-(\hbar e^V_{i j})_{\nu\times\nu}$ is divided into $(m,\ldots,m)$-blocks.
 
\subsection{Notations and Basic Properties}

To explicitly express the (formal) power series solution of system \eqref{Conflu} under the above assumptions, we need to introduce the following notations and results. They can be found in monographs related to analytic ODEs, such as \cite{Balser}, or verified directly.

\begin{nota}
\label{Nota:Mat}
Denote $\tilde{\mathbb{C}}=\{r\mathe^{\mathi \underline{\theta}}: r>0,\theta\in\mathbb{R}\}$ as the universal covering space of $\mathbb{C}\setminus\{0\}$. For a matrix $X$ and an invertible matrix $D$, denote $\tmop{Ad}(D)X\assign D X D^{-1}$.
Denote the operator $A^{\bm{r}}$ and $B^{\bm{l}}$ as 
\begin{align}\label{Def:rlAct}
	A^{\bm{r}}\cdot X \assign X A,\quad
	X\cdot B^{\bm{l}} \assign B X.
\end{align}
For a matrix-valued function $F(z)=\sum_{p=0}^\infty F_pz^p$,
denote
\begin{align}\label{Def:rlFunAct}
	F(A^{\bm{r}})\cdot X \assign \sum_{p=0}^\infty F_p X A^p,\quad
	X\cdot F(B^{\bm{l}}) \assign \sum_{p=0}^\infty B^p X F_p.
\end{align}
For $\nu\times \nu$ (block) matrix $A=(A_{ij})_{\nu\times \nu}$,
\begin{itemize}
\item denote $A_{\ast j}$ as the $j$-th column of $A$;
\item denote $A_{\hat{i}j}$ as the $j$-th column of $A$
      but delete the $i$-th row;
\item denote $A_{\hat{i}\hat{j}}$ as the matrix $A$
      with the $i$-th row and $j$-th column;
\item denote $A_{\hat{i}\ast}$ as the matrix $A$
      with the $i$-th row,
\end{itemize}
and similarly define $A_{i \ast}, A_{i \hat{j}}, A_{\ast\hat{j}}$.
For example, we have
\[ A = \left(\begin{array}{cc}
     A_{11} & A_{1 \hat{1}}\\
     A_{\hat{1} 1} & A_{\hat{1} \hat{1}}
   \end{array}\right) = \left(\begin{array}{cc}
     A_{\hat{\nu} \hat{\nu}} & A_{\hat{\nu} \nu}\\
     A_{\nu \hat{\nu}} & A_{\nu \nu}
   \end{array}\right). \]
Denote the projection operator $\delta_u$ on $\nu\times \nu$ block matrix $A$ as
\begin{align}
\delta_u A & \assign \diag(A_{11},\ldots,A_{\nu\nu}).
\end{align}
\end{nota}


\begin{defi}
Let $\mathrm{Eigen}(A)$ denote the set of
eigenvalues of matrix $A$. If
\begin{align*}
\left(
\mathrm{Eigen}(A)-\mathrm{Eigen}(A)
\right) \cap \mathbb{Z}
= \{0\},
\end{align*}
i.e. the difference between two eigenvalues of $A$
does not take non-zero integers,
we call $A$ \textbf{non-resonant}.
Otherwise, $A$ is called \textbf{resonant}.
\end{defi}

\begin{defi}\label{Def:LkResM}
Define the $k$-th \textbf{recursive matrix} $L_k(z)$ of system \eqref{Conflu} for $k=0,1,\ldots,\nu$, as a 
$(\nu-1)\times(\nu-1)$-block matrix with indices 
$(1,\cdots,k-1,k+1,\cdots,\nu)$,
\begin{align}
\label{LkDef}
L_k (z) = 
\left\{\begin{array}{ll}
     \frac{1}{u_k I - u_{\hat{k} \hat{k}}} 
     \left( (z I + A_{\hat{k} \hat{k}})
     - A_{\hat{k} k} \frac{1}{z I + A_{k k}} A_{k \hat{k}} \right) & ; k = 1,
     \ldots, \nu\\
     - u^{- 1} (z I + A) & ; k = 0
   \end{array}\right..
\end{align}
\end{defi}


The following proposition is then given by a direct computation.
\begin{pro}[Formal solutions at irregular singularity]
\label{FormSolirr}
If $A_{11},\ldots,A_{\nu\nu}$ are all non-resonant, then system \eqref{Conflu} has a unique formal series solution of the following form
\begin{align}
{F}^{[\infty]}(z) = {H}^{[\infty]}(z) \cdot z^{\delta_u A} \mathe^{u z},
\end{align}
where 
\begin{equation*}
{H}^{[\infty]}(z) = 
I+ \sum_{p=1}^\infty {H}_p^{[\infty]} z^{-p},
\end{equation*}
is a formal power series of $n\times n$ matrices. 
The solution ${F}^{[\infty]}(z)$ is
called the \textbf{canonical formal fundamental solution at $z=\infty$}.
For $k=1,\ldots,\nu$, we have
\begin{subequations}
\begin{align}
\label{Cof:Finfhatkk}
({H}_p^{[\infty]})_{\hat{k}k} & =
\left(
\overleftarrow{\prod_{m=1}^{p-1}}
L_k(m-A_{kk}^{\bm{r}})
\right)\cdot
\left(
\frac{1}{u_kI-u_{\hat{k}\hat{k}}}A_{\hat{k}k}
\right),\\
\label{Cof:Finfkk}
({H}_p^{[\infty]})_{kk} & = -
\frac{1}{p-A_{kk}^{\bm{r}}+A_{kk}} A_{k\hat{k}}
\left(
\overleftarrow{\prod_{m=1}^{p-1}}
L_k(m-A_{kk}^{\bm{r}})
\right)\cdot
\left(
\frac{1}{u_kI-u_{\hat{k}\hat{k}}}A_{\hat{k}k}
\right),
\end{align}
\end{subequations}
and we called $\delta_u A$ the \textbf{formal monodromy matrix} of ${F}^{[\infty]}(z)$, where $L_k(z)$ is defined in \eqref{LkDef},
$L_k(m-A_{kk}^{\bm{r}})$ is defined in \eqref{Def:rlFunAct}.
\end{pro}

\begin{defi}
Define the \textbf{anti-Stokes lines} of system \eqref{Conflu} as the rays for which the arguments take the following values
\begin{align}
\label{antiStokes}
\tmop{aS}(u)\assign
\underset{i\neq j}{\bigcup_{1\leqslant i,j\leqslant \nu}}
\left(-\arg(u_i-u_j)+2\pi\mathbb{Z}
\right)
\subseteq\mathbb{R},
\end{align}
The elements of \eqref{antiStokes} are denoted as $(\tau_i)_{i\in\mathbb{Z}}$ and are arranged as $\cdots< \tau_i < \tau_{i+1}<\cdots$,
and are called the \textbf{anti-Stokes arguments}.
\end{defi}
\noindent
It is direct to see that
the anti-Stokes arguments $(\tau_i)_{i\in\mathbb{Z}}$
are periodic. 
That is, there exists an $l$ 
such that for any index $i$,
we have $\tau_{i+l}=\tau_i+\pi$.

\begin{pro}[Analytic solutions at irregular singularity, \cite{Balser}]
\label{Pro:CAFS}
For $d\in(\tau_i,\tau_{i+1})$, there exists a unique analytic function ${H}_d(z)$ on $\tilde{\mathbb{C}}$, 
which is asymptotic to ${H}^{[\infty]}(z)$ on the sector
\begin{align}
\tmop{Sect}_d \assign 
\left\{z\in\tilde{\mathbb{C}}: \arg z \in \left(\tau_i - \frac{\pi}{2}, \tau_{i+1} + \frac{\pi}{2}\right)
\right\},
\end{align}
such that the analytic function
\begin{eqnarray}
{F}_d(z)\assign {H}_d(z) \cdot z^{\delta_u A}\mathe^{u z},
\end{eqnarray}
is a solution of system \eqref{Conflu}.
We call ${F}_d(z)$ 
the \textbf{canonical (analytic) fundamental solution at $z=\infty$}.
\end{pro}

The first part of the following proposition is well known from the general theory of linear differential equation (see e.g., \cite[Chapter 2]{Balser}) and the second part is given by a direct computation.
\begin{pro}[Solutions at regular singularity]
\label{FormSolreg}
If $A$ is non-resonant, then system \eqref{Conflu} has a unique formal series solution of the following form
\begin{align}\label{Series:F0}
{F}^{[0]} (z) = 
{H}^{[0]} (z) \cdot z^A =
\left(
I+ \sum_{p=1}^\infty {H}_p^{[0]} z^{p}
\right) \cdot z^A,
\end{align}
The series \eqref{Series:F0} is a convergent series with an infinite radius of convergence, called the \textbf{canonical fundamental solution at $z=0$}. If we introduce
\begin{align}
{H}^{[k]}(z) & = I+ \sum_{p=1}^\infty {H}_p^{[k]} z^{p}
\assign {H}^{[0]}(z) \mathe^{-u_k z},\quad
k = 0,1, \ldots, \nu,
\end{align}
then we have for $k=1,\ldots,\nu$, 
\begin{subequations}
\begin{align}
\label{Cof:Fkhatkast}
({H}^{[k]}_p)_{\hat{k}\ast} &=
\left(
\overleftarrow{\prod_{m=-p}^{-1}}
L_k(m-A^{\bm{r}})
\right)^{-1}
\cdot I_{\hat{k}\ast},\\
\label{Cof:Fkkast}
({H}^{[k]}_p)_{k\ast} &= -
\frac{1}{-p-A^{\bm{r}}+A_{kk}} A_{k\hat{k}}
\left(
\overleftarrow{\prod_{m=-p}^{-1}}
L_k(m-A^{\bm{r}})
\right)^{-1}
\cdot I_{\hat{k}\ast};
\end{align}
for $k=0$, we conveniently assume that $u_k=0$, and we have
\begin{align}
\label{Cof:F0}
{H}_p^{[0]} = 
\left(
\overleftarrow{\prod_{m=-p}^{-1}}
L_0(m-A^{\bm{r}})
\right)^{-1}
\cdot I,
\end{align}
\end{subequations}
where $L_0(z)$ is defined in \eqref{LkDef}.
\end{pro}

For convenience, we shall always set
\begin{align}
	{H}^{[k]}_0=I,\quad k=0,1,\ldots,\nu,\infty.
\end{align}

\begin{defi}(Monodromy data)\label{monodata}
Suppose that $A_{11},\ldots,A_{\nu\nu}$ are all non-resonant, $d\in\mathbb{R}\setminus{\rm aS}(u)$ do not take the anti-Stokes arguments of system \eqref{Conflu},
$\tau\in {\rm aS}(u)$ take the anti-Stokes arguments.
\begin{itemize}
	\item Denote the following $\nu\times \nu$ block constant matrices
	\begin{subequations}
		\begin{align}
		S_{[\tau]}(u,A) & \assign 
        {F}_{\tau+\varepsilon}(z)^{-1} 
        {F}_{\tau-\varepsilon}(z),\quad
        \tau\in{\rm aS}(u),\quad
        \varepsilon>0
        \text{ sufficiently small},\\
		S_d^\pm(u,A) & \assign {F}_{d\pm\pi}(z)^{-1} {F}_{d}(z),\quad
        d\notin{\rm aS}(u),\\
		S_d(u,A) & \assign S_d^+(u,A)-S_d^-(u,A).
        \label{Def:Sd}
	    \end{align}
	\end{subequations}
	We call $S_{[\tau]}$ the \textbf{Stokes matrix} 
    with respect to the anti-Stokes direction $\tau$,
    and $S_d$ the \textbf{(normalized) Stokes matrix} (in direction $d$);
	\item We additionally assume that $A$ is non-resonant. Denote the following $\nu\times \nu$ block constant matrices
	\begin{align}
		{C}_d(u,A)&\assign {F}_d(z)^{-1}{F}^{[0]}(z),\\
		\label{Def:Md}
		\mathe^{2\pi\mathi M_d(u,A)}&\assign {F}_d(z)^{-1} {F}_d(z\mathe^{2\pi\mathi}),
	\end{align}
	where $M_d$ is the unique matrix that satisfies \eqref{Def:Md} and $\mathrm{Eigen}(M_d)=\mathrm{Eigen}(A)$. We call ${C}_d$ the \textbf{central connection matrix} (in direction $d$), $M_d$ the \textbf{monodromy matrix} (in direction $d$), and $\mathe^{2\pi\mathi M_d}$ the \textbf{monodromy factor} (in direction $d$).
\end{itemize}
When there is no ambiguity, we will omit $u,A$ for $S_{[\tau]}(u,A), S_d^\pm(u,A),C_d(u,A),M_d(u,A)$.
\end{defi}

The Stokes matrices
of the system 
at the irregular singularity fully reflects
the Stokes phenomenon of its canonical fundamental solutions.
This paper will frequently use the following properties of these monodromy data, which can be directly verified from the definition.

\begin{pro}\label{Prop:MD}
Suppose that $A_{11},\ldots,A_{\nu\nu}$ are all non-resonant, $d\in\mathbb{R}$ does not coincide with the anti-Stokes arguments of system \eqref{Conflu}, then we have
\begin{subequations}
\begin{align}
  \label{LUDecomp}
    \mathe^{2\pi\mathi M_d}
    & = 
    (S^-_d)^{- 1} \cdot \mathe^{2 \pi \mathi
    \delta_u A} \cdot S^+_d, \\
    \label{StokesdM}
    S^{\pm}_{d + 2 k \pi}
    & = 
    \mathe^{- 2 k \pi \mathi \delta_u
    A} \cdot S^{\pm}_d \cdot \mathe^{2 k \pi \mathi \delta_u
    A}, \quad k \in \mathbb{Z},\\
    \label{Spm}
    S^{\pm}_{d \mp \pi} 
    & = 
    (S^{\mp}_d)^{- 1},\\
    \label{ShiftStokes}
    S_d^\pm(c u+c_0I,\tmop{Ad}(D)A+c_1I) & =
    \tmop{Ad}(Dc^{\delta_uA})
    S_{d+\arg c}^\pm(u,A),\quad
    c\neq 0,
    \delta_uD=D,\\
    \label{TStokes}
    S_d^\pm(-u^\top,-A^\top) & =
    S_d^\pm(u,A)^{-\top}.
\end{align}
\end{subequations}
If we additionally assume that $A$ is non-resonant, then 
\begin{subequations}
\begin{align}
\label{Md=ConjPhi}
M_d
& = 
{C}_d \cdot A \cdot 
{C}_d^{- 1}, \\
{C}_{d + 2 k \pi}
& = 
\mathe^{-2k\pi\mathi\delta_u A} \cdot 
{C}_d\cdot
\mathe^{2k\pi\mathi A}, \quad k
\in \mathbb{Z}, \\
\label{Cpm}
{C}_{d \pm \pi}
& = 
S^{\pm}_d\cdot {C}_d,\\
    \label{ShiftOmega}
    C_d(c u+c_0I,\tmop{Ad}(D)A+c_1I) & =
    Dc^{\delta_uA}\cdot
    C_{d+\arg c}(u,A)\cdot
    c^{-A}D^{-1},\quad
    c\neq 0,
    \delta_uD=D,\\
    \label{TOmega}
    C_d(-u^\top,-A^\top) & =
    C_d(u,A)^{-\top}.
\end{align}
\end{subequations}
\end{pro}

\begin{pro}[see e.g., \cite{Balser}]\label{AppPro:Sto[t]}
For $j\neq k$ and
anti-Stokes arguments $(\tau_i)_{i\in\mathbb{Z}}$,
\begin{itemize}
    \item the diagonal blocks of $S_{[\tau_i]}$ are identity matrices;
    \item if $-\arg(u_k-u_j)\notin \tau_i+2\pi\mathbb{Z}$, 
          then $(S_{[\tau_i]})_{jk}=0$;
    \item if $-\arg(u_k-u_j)\in \tau_i+2\pi\mathbb{Z}$, 
          then $(S_{[\tau_i]})_{jk}=(S_{\tau_i\pm\varepsilon})_{jk}$,
          where $\varepsilon>0$ sufficiently small;
    \item take $l\in\mathbb{Z}$ such that $\tau_{i+l}=\tau_i+\pi$.
          If $d\in(\tau_i,\tau_{i+1})$, then
          \begin{align}
              S_d^+ = S_{[\tau_{i+l}]} \cdots
                      S_{[\tau_{i+2}]}S_{[\tau_{i+1}]},\quad
              S_d^- = S_{[\tau_{i-l+1}]}^{-1} \cdots
                      S_{[\tau_{i-1}]}^{-1}S_{[\tau_{i}]}^{-1}.
          \end{align}
\end{itemize}
\end{pro}

\begin{pro}[see e.g., \cite{Balser}]\label{AppPro:StoTri}
Suppose that $d\in\mathbb{R}$ does not take the anti-Stokes arguments of system \eqref{Conflu}, and
    \begin{equation}
    \label{Cond:+UpperTri}
      \tmop{Im} (u_{1} \mathe^{\mathi d}) >
      \tmop{Im} (u_{2} \mathe^{\mathi d}) > \cdots > 
      \tmop{Im} (u_{\nu}  \mathe^{\mathi d}),
    \end{equation}
    \begin{itemize}
    \item then the Stokes matrices $S_d^+, S_d^-$ are block upper triangular matrix and lower triangular matrix respectively, with the identity matrix as the diagonal block;
    \item we have
    \begin{subequations}
	\begin{align}\label{Sdii+1}
        (S_d)_{i,i+1} & = (S_{[\tau_+]})_{i,i+1},\quad
        \tau_+=-\arg(u_{i+1}-u_i)\in(d,d+\pi),\\
        \label{Sdi+1i}
        (S_d)_{i+1,i} & = (S_{[\tau_-]})_{i+1,i},\quad
        \tau_-=-\arg(u_i-u_{i+1})\in(d-\pi,d).
    \end{align}
    \end{subequations}
    \end{itemize}
\end{pro}
Proposition \ref{AppPro:Sto[t]} and Proposition \ref{AppPro:StoTri} each provide the basic structure and fundamental relations of $S_{[t]}$ and $S_d$, respectively.
For example, in the generic case, each $S_d^\pm$ has $\frac{n(n-1)}{2}$ nonzero entries. Therefore, the Stokes matrix $S_d$ defined by \eqref{Def:Sd} generally has $n^2-n$ nonzero entries other than diagonal entries, and it can directly recover $S_d^\pm$.
In Section \ref{Sect:PathAlg}, 
we will introduce how to recover all the Stokes matrices $S_{[t]}$ 
from $S_d$ for a given $d$.
If condition \eqref{Cond:+UpperTri} is not satisfied for a direction $d\in\tmop{aS}(u)$, then there exists a unique permutation of matrix indices such that \eqref{Cond:+UpperTri} is valid
and the corresponding Stokes matrices $S_d^+, S_d^-$ are (blocked) upper and lower triangular matrices.

\subsection{Borel-Laplace Transform}

The definitions of the Borel transform and the Laplace transform used in this paper will be slightly different from the ones commonly used.
Under our definitions, the Borel transform and the Laplace transform of the solution to the system \eqref{Conflu} both satisfy the same new system \eqref{RegSys}.

\begin{defi}\label{Def:BL}
	Define the \textbf{formal Borel transform} $\mathcal B$ at $z=\infty$ and the \textbf{formal Laplace transform} $\mathcal L$ at $z=0$ by 
\begin{align}
\label{Def:FBd}
\mathcal B \left(
\sum_{p=0}^\infty a_p z^{-p-s}\mathe^{u_0 z}
\right)
& \assign
\sum_{p=0}^\infty a_p \frac{(u_0-\xi)^{p+s}}{(p+s)!},
\quad u_0\in\mathbb{C},\\
\label{Def:FLd}
\mathcal L \left(
\sum_{p=0}^\infty a_p z^{p-s}\mathe^{u_0 z}
\right)
& \assign
\sum_{p=0}^\infty a_p (\xi-u_0)^{-p+s}(p-s-1)!,
\quad u_0\in\mathbb{C}.
\end{align}
For the analytic function $f(z)$ on the logarithmic Riemann surface around $z = \infty$, define the \textbf{analytic Borel transform} $\mathcal B_d$ with respect to direction $d$ by
\begin{align}
\label{Def:ABd}
(\mathcal B_d f)(\xi)\assign
\frac{1}{2\pi\mathi}\left(
\int
_{N\mathe^{\mathi(d+\frac{\pi+\varepsilon}{2})}}
^{\infty\mathe^{\mathi(d+\frac{\pi+\varepsilon}{2})}} +
\int_\gamma +
\int
^{N\mathe^{\mathi(d-\frac{\pi+\varepsilon}{2})}}
_{\infty\mathe^{\mathi(d-\frac{\pi+\varepsilon}{2})}}
\right)
f(z) \mathe^{-\xi z}\frac{\mathd z}{z},
\end{align}
where $N>0$ is sufficiently large, $\varepsilon>0$ is sufficiently small. The integration path is taken as follows
\[ 
   \raisebox{-0.00217500276123397\height}{\includegraphics[width=7.20864489046307cm,height=3.8600616555162cm]{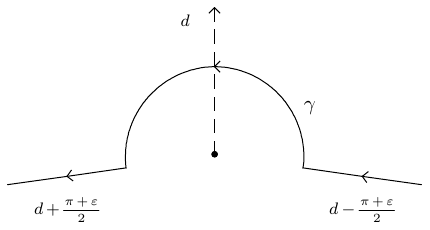}}
\]
For the analytic function $f(z)$ 
on the logarithmic Riemann surface around $z = 0$, define the \textbf{analytic Laplace transform} $\mathcal L_d$ with respect to direction $d$ by
\begin{align}
\label{Def:ALd}
(\mathcal L_d f)(\xi)\assign
\int_{0}^{\infty\mathe^{\mathi d}}
f(z) \mathe^{-\xi z}\frac{\mathd z}{z}.
\end{align}
\end{defi}

It can be seen that for an appropriate complex number $\xi$, both \eqref{Def:ABd} and \eqref{Def:ALd} are meaningful without the need to impose any argument requirements on $u_0-\xi$ or $\xi-u_0$. However, to make the formal transformation \eqref{Def:FBd} and \eqref{Def:FLd} valid, $\xi$ cannot be just some complex numbers, we must specify the arguments of $u_0-\xi$ and $\xi-u_0$, respectively. Therefore, in order to establish the connection between the formal Borel/Laplace transform and the analytic Borel/Laplace transform, it is necessary to clarify the argument conditions as follows.

\begin{lem}
\label{Lem:Borel}
Suppose that $\tilde{f}(z) = \sum_{p=0}^\infty a_p z^{-p-s}\mathe^{u_0 z}$ is a formal series at $z=\infty$ such that $(\mathcal B \tilde{f})(\xi)$ is a convergent series at $\xi= u_0$, and an analytic function $f(z)\mathe^{-u_0 z}$ is
asymptotic to $\tilde{f}(z)\mathe^{-u_0 z}$ on the sectorial region $\{z\in\tilde{\mathbb{C}}: |z|>N, \arg z \in \left(\tau_i - \frac{\pi}{2}, \tau_{i+1} + \frac{\pi}{2}\right)\}$. If $d\in(\tau_i,\tau_{i+1})$, then
\begin{itemize}
	\item $(\mathcal B_d f)(\xi)$ is a convergent integral for $\arg (u_0-\xi)\in (-d-\frac{\varepsilon}{2},-d+\frac{\varepsilon}{2})$ with sufficiently small $\varepsilon$, and is 
    compatible with each other when $d$ varies;
	\item by assigning the argument $\arg (u_0-\xi)\in (-\tau_{i+1},-\tau_i)$ in series $(\mathcal B \tilde{f})(\xi)$ for sufficiently small $u_0-\xi$, we have
\begin{align}
(\mathcal B \tilde{f})(\xi) = (\mathcal B_d f)(\xi).
\end{align} 
\end{itemize}
\end{lem}

\begin{lem}
\label{Lem:Laplace}
Suppose that $f(z)=\sum_{p=0}^\infty a_p z^{p-s}\mathe^{u_0 z}$ is a (formal) series at $z=0$ such that $(\mathcal L f)(\xi)$ is a convergent series at $\xi=\infty$, and $\tmop{Re}s<1$. Then
\begin{itemize}
	\item $(\mathcal L_d f)(\xi)$ is a convergent integral for $\arg \xi\in (-d-\frac{\varepsilon}{2},-d+\frac{\varepsilon}{2})$ with sufficiently large $\xi$, sufficiently small $\varepsilon$, and is 
    compatible with each other when $d$ varies;
	\item by assigning the argument $\arg (\xi-u_0)\in (-d-\frac{\varepsilon}{2},-d+\frac{\varepsilon}{2})$ in series $(\mathcal L f)(\xi)$ with sufficiently large $\xi$, we have
\begin{align}
(\mathcal L f)(\xi) = (\mathcal L_d f)(\xi).
\end{align} 
\end{itemize}
\end{lem}
One checks that the formal/analytic solutions ${F}(z)$ of the confluent hypergeometric system \eqref{Conflu}, after the formal/analytic Borel/Laplace transform, will be transformed into the formal/analytic solutions $Y(\xi)$ of the system
\begin{align}\label{RegSys}
(u-\xi I) \frac{\mathd}{\mathd \xi} Y(\xi) = A \cdot Y(\xi).
\end{align}
Therefore, the analytic solution $Y(\xi)$ is defined on the universal covering space of $\mathbb{C}\setminus\{u_1,\ldots,u_\nu\}$, denoted as $\tilde{\mathbb{C}}(u)$.
We will parameterize the elements in $\tilde{\mathbb{C}}(u)$ in two equivalent ways. 
One is by considering the arguments $\arg(u_k-\xi)$. The space $\tilde{\mathbb{C}}(u)$ with such parametrization is denoted by $\tilde{\mathbb{C}}^{\mathcal B}(u)$; the other is by considering the arguments $\arg(\xi-u_k)$. The same space $\tilde{\mathbb{C}}(u)$ with the latter parametrization is denoted by $\tilde{\mathbb{C}}^{\mathcal L}(u)$.
Furthermore, we will consider the following principal branch
\begin{align}
\mathbb{C}_d(u) \assign \mathbb{C}\setminus 
\bigcup_{k=1}^\nu \{\xi\in\mathbb{C}: \arg(\xi-u_k)=-d\},
\end{align}
and the parameterizations defined on this set (here we use the upper indices to stress the different parameterization on the same space)
\begin{subequations}
\begin{align}
    \mathbb{C}_d^{\mathcal B}(u) &\assign 
    \{\xi\in \tilde{\mathbb{C}}^{\mathcal B}(u):
    \arg(u_k-\xi)\in(-d-\pi,-d+\pi),
    \text{ for all }k=1,\ldots,\nu\},\\
    \mathbb{C}_d^{\pm}(u) &\assign 
    \{\xi\in \tilde{\mathbb{C}}^{\mathcal L}(u):
    \arg(\xi-u_k)\in(-(d\pm\pi)-\pi,-(d\pm\pi)+\pi),
    \text{ for all }k=1,\ldots,\nu\}.
\end{align}
\end{subequations}
By applying Lemma \ref{Lem:Borel} and \ref{Lem:Laplace},
we obtain the following result:
\begin{itemize}
	\item If we parameterize $\mathbb{C}_d(u)$ by 
          $\mathbb{C}_d^{\mathcal B}(u)$, 
          i.e., we assign the argument 
          $\arg(u_k-\xi)\in(-d-\frac{\varepsilon}{2},-d+\frac{\varepsilon}{2})$, then
            \begin{subequations}
        	\begin{align}\label{BdXdk}
	           	(\mathcal B_d {F}_d)(\xi)_{\ast k} & = 
	           	\sum_{p=0}^\infty ({H}_p^{[\infty]})_{\ast k}
	           	\frac{(u_k-\xi)^{p I-A_{kk}}}{(p I-A_{kk})!},\\
	        \label{BdXd}
	           	(\mathcal B_d {F}_d)(\xi) & = 
	           	\sum_{p=0}^\infty {H}_p^{[\infty]}
	           	\frac{(u-\xi)^{p I-\delta_u A}}{(p I-\delta_u A)!}.
        	\end{align}
            \end{subequations}	
Lemma \ref{Lem:Borel} ensures that the left side of \eqref{BdXdk} is a convergent definite integral, and the right side of \eqref{BdXdk} is a convergent power series at $\xi = u_k$, which clarifies the meaning of \eqref{BdXdk}. We interpret \eqref{BdXd} as the analytic continuation of \eqref{BdXdk};
	\item If we parameterize $\mathbb{C}_{d\pm\pi}(u)$ by $\mathbb{C}_{d\pm\pi}^{\mp}(u)$, i.e., we assign the argument $\arg(\xi-u_k)\in(-d-\frac{\pi}{2},-d+\frac{\pi}{2})$, then
	\begin{align}
	\label{LdX0}
	(\mathcal L_{d} {F}^{[0]})(\xi)
	& = \sum_{p=0}^\infty {H}_p^{[k]}((p-1)I+A)! (\xi-u_k)^{-(p I+A)},
	\quad k=0,1,\ldots,\nu,\quad
    \tmop{Re} \mathrm{Eigen}(A)>0,
	\end{align}
	where $\xi$ needs to be sufficiently large, and we conveniently assume that $u_0=0$.
\end{itemize}
It is directly to 
see that condition $\tmop{Re} \mathrm{Eigen}(A)>0$ is not essential, 
as it can always be achieved by replacing the residue matrix $A$ with $A+cI$ instead.

\subsection{Darboux's Method and Path Class with Arguments}

The
Darboux's method in this subsection refers to \cite{Balser88}. 
It was referred to the associate function method in \cite[Chapter 9]{Balser}.
In this subsection, 
for canonical fundamental solution $F_d,F^{[0]}$,
we will focus on the asymptotic behavior of $(\mathcal B_d {F}_d)(\xi),(\mathcal L_{d} {F}^{[0]})(\xi)$ at the singularities $\xi=u_1,\ldots,u_\nu$ 
after analytic continuation.
On the one hand, Proposition \ref{Pro:gammajk} states that these asymptotics can be used to give the monodromy data of the system \eqref{Conflu}. On the other hand, according to Darboux's method (Lemma \ref{Lem:AnaLim}), the same asymptotics can be expressed by the coefficients ${H}_p^{[k]}$. Therefore, we will ultimately be able to directly derive those monodromy data through the recursive matrix $L_k$ applying \eqref{Cof:Finfhatkk} and \eqref{Cof:Fkhatkast}.
First, it follows from Definition \ref{Def:BL} that
\begin{lem}\label{LdFdnice}
    $(\mathcal L_{d\pm\frac{\pi+\varepsilon}{2}}{F}_{d\pm\pi})(\xi)_{\ast t}$ 
    is analytic on $\mathbb{C}\setminus\{\xi\in\mathbb{C}: \arg(\xi-u_t)=-d\}$.
\end{lem}

Second, we have
\begin{lem}\label{Lem:BdFdnice}
There exists a function $\tmop{hol}(u_s-\xi)$
holomorphic at $\xi=u_s$ such that
on $\mathbb{C}^{\mathcal B}_d(u)$
    \begin{align}\nonumber
    (\mathcal B_d {F}_d)(\xi)_{st} 
    & = (I+O(u_s-\xi))
    (u_s-\xi)^{-A_{ss}}(A_{ss}-I)!\cdot
    \frac{\mathe^{\pi\mathi A_{ss}} (S_d^+)_{st}-
	\mathe^{-\pi\mathi A_{ss}} (S_d^-)_{st}}{2\pi\mathi}\\
    & \quad
    +\tmop{hol}(u_s-\xi),\quad \xi\to u_s.
    \label{BdXdjkgamma}
    \end{align}
\end{lem}

\begin{prf}
Without loss of generality, we assume that
\begin{align}
    \tmop{Re} \mathrm{Eigen}(A)>0, \quad 
    \mathrm{Eigen}(\delta_u A)\cap \mathbb{Z}=\varnothing.
\end{align}
Condition 
$\tmop{Re} \mathrm{Eigen}(A)>0$ 
ensures that
\begin{align}\nonumber
    \mathcal B_d {F}_d & =
    \frac{1}{2\pi\mathi}\left(
    \mathcal L_{d+\frac{\pi+\varepsilon}{2}} -
    \mathcal L_{d-\frac{\pi+\varepsilon}{2}} \right)
    {F}_{d}
    = \frac{1}{2\pi\mathi}
    \mathcal L_{d\pm\frac{\pi+\varepsilon}{2}} {F}_d \cdot
    \pm(I-\mathe^{\mp 2\pi\mathi M_d}).
\end{align}
Combining \eqref{Md=ConjPhi} in Proposition \ref{Prop:MD}, we have
\begin{align}
    \mathcal B_d {F}_d & =
    \frac{1}{2\pi\mathi}
    \mathcal L_{d\pm\frac{\pi+\varepsilon}{2}} {F}_{d\pm\pi} \cdot
    \mathe^{\mp\pi\mathi\delta_u A}
    (\mathe^{\pi\mathi\delta_u A}S_d^+-
    \mathe^{-\pi\mathi\delta_u A}S_d^-).
\end{align}
From Lemma \ref{LdFdnice} and \eqref{BdXdk}, 
on $\mathbb{C}^{\mathcal B}_d(u)$ we have
\begin{subequations}
\begin{align}
	(\mathcal B_d {F}_d)_{st} & =
	\frac{1}{2\pi\mathi}
	(\mathcal L_{d\pm\frac{\pi+\varepsilon}{2}} {F}_{d\pm\pi})_{ss} \cdot
	\mathe^{\mp\pi\mathi A_{ss}}
	(\mathe^{\pi\mathi A_{ss}} (S_d^+)_{st}-
	\mathe^{-\pi\mathi A_{ss}} (S_d^-)_{st})
	+\tmop{hol}(u_s-\xi),
	\label{BdXdjk}\\
	\nonumber
	(\mathcal B_d {F}_d)_{ss} & =
	\frac{1}{2\pi\mathi}
	(\mathcal L_{d\pm\frac{\pi+\varepsilon}{2}} {F}_{d\pm\pi})_{ss} \cdot
	\mathe^{\mp\pi\mathi A_{ss}}
	(\mathe^{\pi\mathi A_{ss}} -
	\mathe^{-\pi\mathi A_{ss}} )
	+\tmop{hol}(u_s-\xi)\\
	& = (I+O(u_s-\xi))
	\frac{(u_s-\xi)^{-A_{ss}}}{(-A_{ss})!},
	\quad \xi\to u_s,
\label{BdXdjj}
\end{align}
\end{subequations}
Condition $\mathrm{Eigen}(\delta_u A)\cap \mathbb{Z}=\varnothing$ ensures that we can substitute \eqref{BdXdjj} into \eqref{BdXdjk}
, and deduce that
\begin{align}
\nonumber
(\mathcal B_d {F}_d)(\xi)_{st} & =
(I+O(u_s-\xi))
\frac{(u_s-\xi)^{-A_{ss}}}{(-A_{ss})!}\cdot
\frac{1}
{\mathe^{\pi\mathi A_{ss}} - \mathe^{-\pi\mathi A_{ss}}}
\left(\mathe^{ \pi\mathi A_{ss}} (S_d^+)_{st}-
\mathe^{-\pi\mathi A_{ss}} (S_d^-)_{st}
\right)\\
& \quad +\tmop{hol}(u_s-\xi),
\end{align}
thus we finish the proof.
\end{prf}

\begin{cor}\label{Sizeatus}
    If $Y(\xi)$ is a solution of system \eqref{RegSys}, then
    \begin{enumerate}
        \item there exists a unique $n_s\times n$-matrix $C_{s\ast}$ and
        an analytic function $\tmop{hol}(\xi-u_s)$ at $\xi=u_s$
        such that
        \begin{align}\label{Ysastnice}
            Y(\xi)_{s\ast} = (I+O(\xi-u_s))
            (\xi-u_s)^{-A_{ss}} \cdot C_{s\ast}
            + \tmop{hol}(\xi-u_s),\quad
            \xi\to u_s;
        \end{align}

        \item there exists a unique $n\times n$-matrix $C$
        such that
        \begin{align}\label{Ynice}
            Y(\xi) = (I+O(-\xi))
            (-\xi)^{-A} \cdot C,
            \quad \xi\to\infty.
        \end{align}
    \end{enumerate}
\end{cor}

\begin{prf}
    The uniqueness is direct; we will only show the existence.
    Since there exists a constant matrix $C$
    such that
    $Y(\xi) = (\mathcal B_d {F}_d)(\xi)\cdot C$, 
    by Lemma \ref{Lem:BdFdnice}, 
    we can prove \eqref{Ysastnice}.
    Similarly, from \eqref{LdX0}, we can obtain \eqref{Ynice}.
\end{prf}

\begin{defi}
\label{Def:StokesWeight}
For $u_1,\ldots,u_\nu\in\mathbb{C}$, $u_\infty\assign\infty$
and $s,t\in\{1,\ldots,\nu,\infty\}$,
take $\gamma_{st}$ as a path
in $\mathbb{C}\setminus\{u_1,\ldots,u_\nu\}$ from $u_s$ to $u_t$ 
with the following real arguments
\begin{subequations}
\begin{align}
    \theta_s & = 
    \left\{\begin{array}{ll}
    \underset{\xi \in \gamma_{s t}, \xi \rightarrow u_s}{\lim} 
    \arg (\xi - u_s), & s \neq \infty\\
    \underset{\xi \in \gamma_{s t}, \xi \rightarrow u_s}{\lim} 
    \arg (c - \xi),\quad\forall c\in\mathbb{C},
    & s = \infty
   \end{array}\right.,\\
    \theta_t & = 
    \left\{\begin{array}{ll}
    \underset{\xi \in \gamma_{s t}, \xi \rightarrow u_t}{\lim}
    \arg (u_t - \xi), & t \neq \infty\\
    \underset{\xi \in \gamma_{s t}, \xi \rightarrow u_t}{\lim} 
    \arg (\xi - c),\quad\forall c\in\mathbb{C},
    & t = \infty
   \end{array}\right..
\end{align}
\end{subequations}
We call the equivalence classes of paths with real arguments 
$\boldsymbol{\gamma}_{st} =
[\gamma_{st},\theta_s,\theta_t]$
under fixed-endpoint homotopy as the 
\textbf{path class with arguments}.
For distinction, 
we always use boldface $\boldsymbol{\gamma}_{st}$
to represent the path class with arguments
and regular typeface $\gamma_{st}$
to represent its representative element.
\end{defi}

Corollary \ref{Sizeatus} and
the non-resonant condition on $A$ and the diagonal block $A_{ss}$ 
ensure that the following definition is well-defined.

\begin{defi}\label{Def:WeightS}
\begin{subequations}
Denote $n_\infty\assign n_1+\cdots+n_\nu$.
If we continue in the reverse direction along $\boldsymbol{\gamma}_{st}$,
then we will take
$\mathcal S_{u,A}(\gamma_{st},\theta_s,\theta_t)$
as an $n_s\times n_t$-matrix defined by
\begin{align}
\nonumber
\sum_{p=0}^\infty ({H}_p^{[\infty]})_{st}
\frac{(u_t-\xi)^{p I-A_{tt}}}{(p I-A_{tt})!}
& = (I+O(\xi-u_s))
(\xi-u_s)^{-A_{ss}}(A_{ss}-I)!\cdot
\frac{\mathcal S_{u,A}(\gamma_{st},\theta_s,\theta_t)}{2\pi\mathi}\\
& \quad +\tmop{hol}(\xi-u_s),\quad
s\neq\infty,\quad t\neq\infty,
\label{Weight:S}\\
\nonumber
\sum_{p=0}^\infty ({H}_p^{[\infty]})_{\ast t}
\frac{(u_t-\xi)^{p I-A_{tt}}}{(p I-A_{tt})!}
& = (I+O(-\xi))
(-\xi)^{-A}(A-I)!\cdot
\frac{\mathcal S_{u,A}(\gamma_{\infty t},\theta_\infty,\theta_t)}{2\pi\mathi},\\
& \quad
s=\infty,\quad t\neq\infty,
\label{Weight:Sinf}\\
\nonumber
\sum_{p=0}^\infty ({H}_p^{[s]})_{s\ast}
\frac{((p-1)I+A)!(\xi-u_s)^{-(pI+A)}}
{(A-I)!(-A)!}
& = (I+O(\xi-u_s))
(\xi-u_s)^{-A_{ss}}(A_{ss}-I)!\cdot
\frac{\mathcal S_{u,A}(\gamma_{s\infty},\theta_s,\theta_\infty)}
{2\pi\mathi}\\
&\quad 
+\tmop{hol}(\xi-u_s),
\quad s\neq\infty,\quad t=\infty,
\label{Weight:O}\\
\nonumber
\sum_{p=0}^\infty {H}_p^{[0]}
\frac{((p-1)I+A)! \xi^{-(pI+A)}}
{(A-I)!(-A)!}
& = (I+O(-\xi))
(-\xi)^{-A}(A-I)!\cdot
\frac{\mathcal S_{u,A}(\gamma_{\infty\infty},\theta_\infty^{(s)},\theta_\infty^{(t)})}
{2\pi\mathi},\\
&\quad 
s=\infty,\quad t=\infty.
\label{Weight:H0p}
\end{align}
\end{subequations}
\end{defi}

\begin{figure}[h!]
\label{fig:StokesPath}
    \centering
    \includegraphics[scale=0.9]{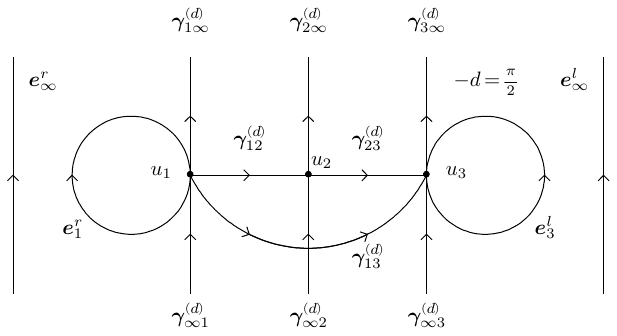}
    \caption{Examples of $\gamma_{st}^{(d)}$
    defined in Proposition \ref{Pro:gammajk}
    for $\nu=3$ 
    \label{Fig:path}}
\end{figure}

\begin{pro}\label{Pro:gammajk}
For distinct $s,t\in\{1,\ldots,\nu\}$,
denote the path class with arguments $\boldsymbol\gamma_{st}^{(d)}$ 
as a path in $\mathbb{C}_d(u)$ from $u_s$ to $u_t$ 
with arguments
\begin{align}\label{Arg:st}
    \theta_s \in 
    \left\{\begin{array}{ll}
     (-d-2\pi,-d), 
     & \tmop{Im} (u_s \mathe^{\mathi d}) > 
     \tmop{Im} (u_t \mathe^{\mathi d})\\
     (-d,-d+2\pi), 
     & \tmop{Im} (u_s \mathe^{\mathi d}) < 
     \tmop{Im} (u_t \mathe^{\mathi d})
    \end{array}\right.,\quad
    \theta_t \in (-d-\pi,-d+\pi),
\end{align}
and denote the path class with arguments 
$\boldsymbol\gamma_{s\infty}^{(d)} \assign 
([u_s,\infty\mathe^{-\mathi d}),-d,-d)$,
$\boldsymbol\gamma_{\infty t}^{(d)} \assign 
((\infty\mathe^{-\mathi d},u_t],-d,-d)$
as a ray (see Figure \ref{Fig:path})
, then we have
\begin{align}\label{BasicWeight}
    \mathcal S_{u,A}(\boldsymbol\gamma_{st}^{(d)})
    = (S_d)_{st},\quad
    \mathcal S_{u,A}(\boldsymbol\gamma_{s\infty}^{(d)})
    = ({C}_d)_{s\ast} \cdot
    (\mathe^{\pi\mathi A} - \mathe^{-\pi\mathi A}),\quad
    \mathcal S_{u,A}(\boldsymbol\gamma_{\infty t}^{(d)})
    = (\mathe^{\pi\mathi A} - \mathe^{-\pi\mathi A})\cdot
    ({C}_d^{-1})_{\ast t}.
\end{align}
Denote the path class with arguments 
($\boldsymbol{e}_s^{l}$) $\boldsymbol{e}_s^r$
as a simple loop
starting from $u_s$ 
in a (counter-) clockwise direction
with arguments $\theta_s=\theta_t$
(see Figure \ref{Fig:path}), 
then we have
\begin{align}\label{LoopWeight}
    \mathcal S_{u,A}(\boldsymbol{e}_s^{l}) 
    = \mathe^{2\pi\mathi A_{ss}} - I,\quad
    \mathcal S_{u,A}(\boldsymbol{e}_s^{r})
    = I - \mathe^{-2\pi\mathi A_{ss}},\quad
    A_{\infty\infty}\assign A.
\end{align}
\end{pro}

\begin{prf}
    From \eqref{BdXdk}, \eqref{LdX0}, we reduce the series in
    \eqref{Weight:S}, \eqref{Weight:O} 
    to $(\mathcal B_d {F}_d)(\xi)$,
    $(\mathcal L_{d} {F}^{[0]})(\xi)
    \frac{1}{(A-I)!(-A)!}$.
    Thus, \eqref{LoopWeight} follows from the definition;
    $\mathcal S_{u,A}(\boldsymbol\gamma_{st}^{(d)})
    = (S_d)_{st}$
    is a direct corollary of Lemma \ref{Lem:BdFdnice}.
    Here, $\xi$ originates from the terminal point $u_t$, 
    and is continued in the reverse direction along $\gamma_{st}^{(d)}$ to obtain the left-hand side of \eqref{BdXdjkgamma}.
    
    We can rewrite \eqref{BdXdjj} as
    \begin{align}
    (\mathcal L_{d\pm\frac{\pi+\varepsilon}{2}} {F}_{d\pm\pi})_{ss}
    & =
    (I+O(\xi-u_s))
    (\xi-u_s)^{-A_{ss}}(A_{ss}-I)!
    +\tmop{hol}(\xi-u_s),
    \end{align}
    where $\arg(\xi-u_s)\in(-(d\pm\pi)-\pi,-(d\pm\pi)+\pi)$.
    Thus we have
    \begin{align}\label{BdX0jastgamma}
    (\mathcal L_{d\pm\frac{\pi+\varepsilon}{2}} {F}^{[0]})_{s\ast}
    & =
    (I+O(\xi-u_s))
    (\xi-u_s)^{-A_{ss}}(A_{ss}-I)!\cdot ({C}_{d\pm\pi})_{s\ast}
    +\tmop{hol}(\xi-u_s).
    \end{align}
    Here, $\xi$ originates from the terminal point $\infty$, and is continued in the reverse direction along $\gamma_{s\infty}^{(d\pm\pi)}$ to obtain the left-hand side of \eqref{BdX0jastgamma}.
    Thus, we have completed the proof of \eqref{BasicWeight}.
\end{prf}

\begin{cor}\label{ProveWeight}
For $\xi\to\infty$,
we have the following connection formula
\begin{align}
    (\mathcal L_{d} {F}^{[0]})(\xi)
    \frac{1}
    {(A-I)!(-A)!} =
    (\mathcal B_d {F}_d)(\xi)\cdot
    C_d\left(\frac{u_k-\xi}{\xi-u_k}
    \right)^A,\quad
    \arg(u_k-\xi)\in(-d-\pi,-d+\pi).
\end{align}
The monodromy factor of $\mathcal B_d {F}_d$ on 
$\mathbb{C}_d^{\mathcal B}(u)$ at $\xi=u_s$ (see Definition \ref{monodata}) is
\begin{align}\label{MonoBdFd}
    U_d^{[m]} = \left\{
    \begin{array}{ll}
    (I-E_{ss}(S_d^\pm -
    \mathe^{\mp 2\pi\mathi\delta_u A} S_d^\mp))^{\pm 1}, 
    & s\neq\infty\\
    \mathe^{-2\pi\mathi M_d}, 
    & s=\infty
    \end{array}\right..
\end{align}
\end{cor}



If we consider a path class with arguments $\boldsymbol{\gamma}$ that can be homotopic to a segment or ray, 
then $\mathcal{S}_{u,A}(\boldsymbol{\gamma})$ can be directly related to the coefficients 
$(H_p^{[\infty]})_{p\in\mathbb{N}}$ or
$(H_p^{[k]})_{p\in\mathbb{N}}$ of the 
corresponding canonical fundamental solutions 
$F^{[\infty]}$ or $F^{[k]}$, respectively.

\begin{lem}[Darboux's method]\cite{Balser88}
\label{Lem:AnaLim}
Let $(c_p)_{p\in\mathbb{N}}$ be a sequence of complex numbers, $\alpha \geqslant 1$ be a real number, and $\beta,C\in\mathbb{C}$.
\begin{itemize}
	\item If the following limit exists
	\begin{align}
	\lim_{p\to\infty}
	p^{-\alpha+1}(v_t-v_s)^{p-\beta}\cdot c_p =
	C\cdot (v_t-v_s)^{-\alpha},
	\end{align}
	then the convergent power series $f(\xi)=\sum_{p=0}^\infty c_p(v_t-\xi)^{p-\beta}$ can be analytically continued along the segment $[v_s,v_t]$ in the reverse direction, such that
	\begin{align}
	\label{Lem:Ana}
	f(\xi)= C\cdot
	(\xi-v_s)^{-\alpha}(\alpha-1)! 
	+ o((\xi-v_s)^{-\alpha}),\quad
	\xi\to v_s,\xi\in[v_s,v_t].
	\end{align}
	Here we fixed the arguments 
    $\theta_s=\theta_t=\arg(v_t-v_s)$,
    making $[v_s,v_t]$ a path class with arguments.
	\item If the following limit exists
	\begin{align}
	\lim_{p\to\infty}
	p^{-\alpha+1}(v_s-v_t)^{-p-\beta}\cdot c_p =
	C\cdot (v_s-v_t)^{-\alpha},
	\end{align}
	then the convergent power series $f(\xi)=\sum_{p=0}^\infty c_p(\xi-v_t)^{-p-\beta}$ can be analytically continued along the ray 
    $[v_s,\infty \mathe^{\mathi \arg(v_s-v_t)})$
    in the reverse direction, such that
	\begin{align}
	f(\xi)= C\cdot
	(\xi-v_s)^{-\alpha}(\alpha-1)!
	+ o((\xi-v_s)^{-\alpha}),\quad
	\xi\to v_s,
	\xi\in [v_s,\infty \mathe^{\mathi \arg(v_s-v_t)}).
	\end{align}
	\[ \raisebox{-0.76491908446748\height}{\includegraphics[width=3.8705562114653cm,height=0.904843893480257cm]{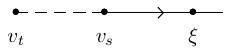}}\]
	Here we fixed the arguments
    $\theta_s=\theta_\infty=\arg(v_s-v_t)$, 
    making $[v_s,\infty \mathe^{\mathi \arg(v_s-v_t)})$ a path class with arguments.
\end{itemize}
\end{lem}
\begin{cor}\label{Weight:Seg}
Suppose that none of the elements $u_1,\ldots,u_\nu$ 
lie inside the segment $[u_s,u_t]$, then we have 
\begin{align}\label{SegWeight}
    \mathcal S_{u,A}([u_s,u_t],-\tau,-\tau) =
    (S_{[\tau]})_{st}.
\end{align}
\end{cor}
\begin{prf}
Since
$[u_s,u_t] \simeq \gamma_{st}^{(\tau\pm\varepsilon)}$,
from \eqref{BasicWeight} in Proposition \ref{Pro:gammajk}
we have $\mathcal S ([u_s,u_t],-\tau,-\tau) = (S_{\tau\pm\varepsilon})_{st}$.
Using Proposition \ref{AppPro:Sto[t]}, we concluded the proof.
\end{prf}

\begin{cor}\label{Cor:LimCoef}
\begin{subequations}
Suppose that $A$, $\delta_u A$ are diagonalizable,
and fix the argument $\arg(u_t-u_s)=-\tau$.
If the following limit exists, 
and none of $u_1,\ldots,u_\nu$
lie inside the segment $[u_s,u_t]$ or 
the ray $[u_t,\infty\mathe^{-\mathi\tau})$,
then we have
\begin{align}
\label{Lim:Stokes}
\lim_{p\to\infty}
(u_t-u_s)^{A_{ss}}p^{-A_{ss}}\cdot
\frac{(u_t-u_s)^p}{(p-1)!}({H}_p^{[\infty]})_{st}\cdot
p^{A_{tt}}(u_t-u_s)^{-A_{tt}} & =
\frac{(S_{[\tau]})_{st}}{2\pi\mathi},
\\\lim_{p\to\infty}
(u_t-u_s)^{A_{tt}}p^{-A_{tt}}\cdot
({H}_p^{[s]})_{t\ast} \frac{p!}{(u_t-u_s)^p}\cdot
p^{A}(u_t-u_s)^{-A} & =
({C}_{\tau\pm\varepsilon})_{t\ast}.
\label{Lim:Omega}
\end{align}
\end{subequations}
\end{cor}
\begin{prf}
Assume that $A_{ss},A_{tt}$ and $A$ have diagonalization
\begin{align*}
P_s A_{ss} P_s^{-1} & =D_s
= \diag(\lambda^{(s)}_1,\ldots,\lambda^{(s)}_{n_s}),\\
P_t A_{tt} P_t^{-1} & =D_t
= \diag(\lambda^{(t)}_1,\ldots,\lambda^{(t)}_{n_s}),\\
P A P^{-1} & =D
= \diag(\lambda_1,\ldots,\lambda_{n}).
\end{align*}
Choose $c^{(s)}_i,c^{(t)}_j\in\mathbb{C}$ such that
$\lambda^{(s)}_i+c^{(s)}_i,
\lambda^{(t)}_j+c^{(t)}_j\geqslant 1$.
By applying Corollary \ref{Weight:Seg}
and replacing the residue matrix $A$ 
with $A+c^{(s)}_iI,A+c^{(t)}_jI$ respectively, 
from \eqref{ShiftStokes} and \eqref{ShiftOmega}
in Proposition \ref{Prop:MD}
we can rewrite \eqref{Weight:S} and \eqref{Weight:O} as
\begin{subequations}
\begin{align}
\nonumber
& \quad \sum_{p=0}^\infty 
(P_s({H}_p^{[\infty]})_{st} P_t^{-1})_{ij} \cdot
\frac{(u_t-\xi)^{p-(\lambda^{(t)}_j+c^{(s)}_i)}}
{(p-(\lambda^{(t)}_j+c^{(s)}_i))!}\\
& =
(\xi-u_s)^{-(\lambda^{(s)}_i+c^{(s)}_i)}((\lambda^{(s)}_i+c^{(s)}_i)-1)!\cdot
\frac{(P_s (S_{[\tau]})_{st} P_t^{-1})_{ij}}{2\pi\mathi}
+O((\xi-u_s)^{1-(\lambda^{(s)}_i+c^{(s)}_i)}),
\label{Weight:SConj}\\
\nonumber
&\quad \sum_{p=0}^\infty 
(P_t({H}_p^{[s]})_{t\ast}P^{-1})_{jl}\cdot
(p-1+(\lambda_l+c^{(t)}_j))!
(\xi-u_s)^{-(p+(\lambda_l+c^{(t)}_j))}\\
& =
(\xi-u_t)^{-(\lambda^{(t)}_j+c^{(t)}_j)}
((\lambda^{(t)}_j+c^{(t)}_j)-1)!\cdot
(P_t ({C}_{\tau\pm\varepsilon})_{t\ast} P^{-1})_{jl}
+O((\xi-u_t)^{1-(\lambda^{(t)}_j+c^{(t)}_j)}),
\label{Weight:OConj}
\end{align}
\end{subequations}
with $\gamma_{st}=[u_s,u_t],\gamma_{t\infty}=[u_t,\infty\mathe^{-\mathi\tau})$.

Let us only assume the existence of the limits in \eqref{Lim:Stokes}, \eqref{Lim:Omega}.
By applying Lemma \ref{Lem:AnaLim} to \eqref{Weight:SConj}, \eqref{Weight:OConj}, 
we can deduce that \eqref{Lim:Stokes}, \eqref{Lim:Omega} holds entrywise, thus completing the proof.
\end{prf}\\
\noindent
We will provide the conditions for the existence of the limit in
Corollary \ref{Cor:LimCoef} in Section \ref{Sect:DiffSys}.

\section{The Associated Difference Systems}\label{Sect:DiffSys}
In this section, we discuss the analytic property of the infinite matrix products in Corollary \ref{Cor:LimCoef}. We provide some explicit convergence conditions of the infinite matrix products, and interpret them as solutions of the closely related $\mu\times\mu$-blocked difference equations for $k=1,...,\nu$
\begin{equation}\label{kdifference}
    \Psi(z+1) = L_k(z) \Psi(z).
\end{equation}

\subsection{Infinite Matrix Product as Solutions of Difference Systems}
Let us first recall a bit more general case. 
Consider the $\mu\times\mu$-blocked difference system
\begin{align}\label{Sys:Difference}
    \Psi(z+1) = v \tilde{B}(z) \Psi(z),
\end{align}
where $v$ is a diagonal matrix, and the coefficient matrix $\tilde{B}(z)$
has a convergent expansion of the following form at $z=\infty$,
\begin{align}\label{Coef:Difference}
    \tilde{B}(z) = z \left(I + \frac{B}{z} +
    \sum_{p=2}^\infty \frac{B_p}{z^p}\right),\quad |z|>R.
\end{align}
Let us further impose the following conditions:
\begin{itemize}
    \item $\tilde{B}(z)$ is a rational function with poles
    $\Lambda = \{\lambda_1,\ldots,\lambda_m\}$
    other than $\infty$;
    \item $v=\diag(v_1I_{n_1},\ldots,v_\mu I_{n_\mu})$ is invertible,
    where $v_1,\ldots,v_\mu$ are distinct, 
    with multiplicity $n_1,\ldots,\nu_\mu$;
    \item Coefficients $B,(B_p)_{p\geqslant2}$ of the system \eqref{Coef:Difference} is divided into
    $(n_1,\ldots,\nu_\mu)$-blocks,
    and we denote $B=(B_{i j})_{\mu\times\mu}$
    as a $\mu\times\mu$-block matrix.
\end{itemize}

\begin{pro}\label{DiffLim}
Suppose that
\begin{align}
    \frac{\varepsilon_B}{2}\assign
    \max\{\|B_{11}\|,\ldots,\|B_{\mu\mu}\|\} < \frac{1}{4}.
\end{align}
For $z\notin \mathbb{Z}^{<0} + \Lambda$, 
if $|v_s|<\varepsilon |v_k|$
for every $k\neq s$, and for some
sufficiently small $\varepsilon(z,\tilde{B})>0$,
then the following limit exists
\begin{align}\label{Lim:DiffPos}
    \lim_{p\to\infty}
	p^{-zI-B_{ss}}
	\frac{v_s^{-p}}{p!}
	\left(
	\overleftarrow{\prod_{m=1}^{p}}
	v \tilde{B}(m+z)
	\right)_{s\ast}.
\end{align}
For $z\notin \mathbb{Z}^{>0} + \Lambda$, 
if $|v_s|<\varepsilon |v_k|$
for every $k\neq s$, and 
sufficiently small $\varepsilon(z,\tilde{B})>0$,
then the following limit exists
\begin{align}\label{Lim:DiffNeg}
    \lim_{p\to\infty}
	\left(
	\overleftarrow{\prod_{m=-p}^{-1}}
	v \tilde{B}(m+z)
	\right)_{\ast s}
	\frac{(-v_s)^{-p}}{p!}
  	p^{zI+B_{ss}}.
\end{align}
\end{pro}

\begin{prf}
    Denote $\mathcal{P}_{st}$ as the set of index sequences 
    $(\mathcal{p}_i)_{i\in\mathbb{N}}\in\{1,\ldots,\mu\}^{\mathbb{N}}$, where
    \begin{align}
        \lim_{i\to\infty} \mathcal{p}_i = s,\quad \mathcal{p}_0=t.
    \end{align}
    For $\mathcal{p}\in\mathcal{P}_{st}$ and $m\in\mathbb{Z}^+$,
    denote 
    \begin{align*}
        \beta^{\mathcal{p}}_m(z) & \assign
        \left( \prod_{k=1}^m \left(
        I+\frac{B_{\mathcal{p}_m\mathcal{p}_m}}{k+z}
        \right)^{-1} \right)
        \left(\delta_{\mathcal{p}_m\mathcal{p}_{m-1}} +
        \frac{B_{\mathcal{p}_m\mathcal{p}_{m-1}}}{m+z} +
        \sum_{p=2}^\infty \frac{(B_p)_{\mathcal{p}_m\mathcal{p}_{m-1}}}{(m+z)^p}
        \right)
        \prod_{k=1}^{m-1} \left(
        I+\frac{B_{\mathcal{p}_{m-1}\mathcal{p}_{m-1}}}{k+z}
        \right),\\
        \beta^{\mathcal{p}}(z) & \assign
        \overleftarrow{\prod_{m=1}^{\infty}} \beta^{\mathcal{p}}_m(z).
    \end{align*}
    Note that for any $m \in \mathbb{N}^+$
    there exist constants $C_0,C_1$ independent of $\mathcal{p}$ such that
    \begin{align}
        \left\lVert
        \prod_{k=1}^m \left(
        I+\frac{B_{\mathcal{p}_m\mathcal{p}_m}}{k+z}
        \right)^{-1}
        \right\rVert & \leqslant
        C_0(z,B) \cdot m^{\lVert B_{\mathcal{p}_m\mathcal{p}_m} \rVert},\\
        \left\lVert
        \prod_{k=1}^{m-1} \left(
        I+\frac{B_{\mathcal{p}_{m-1}\mathcal{p}_{m-1}}}{k+z}
        \right)
        \right\rVert & \leqslant
        C_0(z,B) \cdot m^{\lVert B_{\mathcal{p}_{m-1}\mathcal{p}_{m-1}} \rVert},
    \end{align}
    thus we have
    \begin{align}
        \| \beta^{\mathcal{p}}_m (z) \| & \leqslant 
        \left\{
        \begin{array}{ll}
        \frac{C_1(z, \tilde{B}) }{m^{1 - \varepsilon_B}}, & \mathcal{p}_m \neq
        \mathcal{p}_{m - 1}\\
        1 + \frac{C_1(z, \tilde{B}) }{m^{2 - \varepsilon_B}}, & \mathcal{p}_m =
        \mathcal{p}_{m - 1}
        \end{array}\right.,
    \end{align}
    and $\beta^{\mathcal{p}}(z)$ converges. 
    To prove that the limit \eqref{Lim:DiffPos} exists, 
    we only need to show that the series
    \begin{align}
    \sum_{\mathcal{p}\in\mathcal{P}_{st}}
    \varepsilon^{w(\mathcal{p})}\beta^{\mathcal{p}}(z) =
        \sum_{n=0}^\infty \varepsilon^n
        \sum_{\mathcal{p}\in\mathcal{P}_{st},w(\mathcal{p})=n}
        \beta^{\mathcal{p}}(z)
    \end{align}
    is absolutely convergent,
    where $w(\mathcal{p}) \assign \#\{i\in\mathbb{N^+}:\mathcal{p}_i\neq s\}$.
    Since 
    \begin{align}
        \sum_{\mathcal{p}\in\mathcal{P}_{st},w(\mathcal{p})=n}
        \left\|\beta^{\mathcal{p}}(z)\right\| \leqslant
        2^{n-1}(\mu-1)^n\cdot
        \left(\sum_{m=1}^\infty
        \frac{1}{m^{2 - 2\varepsilon_B}}\right)^{2n}
        C_1(z, \tilde{B} )^{2n}
        C_2(z, \tilde{B} ),
    \end{align}
    the existence of the limit \eqref{Lim:DiffPos} is thus proved.
    The existence of the limit \eqref{Lim:DiffNeg} is similar.
\end{prf}\\
If the number
of distinct diagonal elements of $v$
is $\mu=1$,
then the control condition for $v_s$
in Proposition \ref{DiffLim} will automatically hold.

It is direct to verify that
the limits \eqref{Lim:DiffPos} and \eqref{Lim:DiffNeg} provide
the canonical analytic solutions
$\mathscr{L}^+(z;v_1,\ldots,v_\mu)$,
$\mathscr{L}^-(z;v_1,\ldots,v_\mu)$
to the difference equation \eqref{Sys:Difference}.
\begin{defi}\label{soldifference}
For a fixed invertible $v=\diag(v_1I_{n_1},\ldots,v_\mu I_{n_\mu})$,
we introduce the canonical analytic solutions of \eqref{Sys:Difference},
when the following limits exists
\begin{align}
    \mathscr{L}^+(1+z)^{-1} & \assign
    v^{-(1+z)} \left( \lim_{p\to\infty}
	p^{-zI-\delta_vB}
	\frac{v^{-p}}{p!}
	\left( \overleftarrow{\prod_{m=1}^{p}}
	v \tilde{B}(m+z) \right) \right),\quad
    z \notin \mathbb{Z}^{<0}+\Lambda, \\
    \mathscr{L}^-(z) & \assign
    \left( \lim_{p\to\infty}
	\left( \overleftarrow{\prod_{m=-p}^{-1}}
	v \tilde{B}(m+z)
	\right)
	\frac{(-v)^{-p}}{p!}
  	p^{zI+\delta_vB} \right) (-v)^z,\quad
    z \notin \mathbb{Z}^{>0}+\Lambda.
\end{align}
When the convergence conditions are invalid, the corresponding analytic functions obtained through analytic continuation are also denoted as
$\mathscr{L}^+(1+z)^{-1}$,
$\mathscr{L}^-(z)$.
\end{defi}

\begin{rmk}
The canonical solutions given in Definition \ref{soldifference}
are the same as those discussed by Birkhoff \cite{birkhoff1911general}, 
although the forms of the definitions differ slightly.
Therefore, according to \cite{birkhoff1911general}, 
the only singularities of $\mathscr{L}^+(z)^{-1}$ and $\mathscr{L}^-(z)$ 
are precisely 
$\mathbb{Z}^{\leqslant 0}+\Lambda$ and 
$\mathbb{Z}^{>0}+\Lambda$.
\end{rmk}

Now back to our equation \eqref{kdifference}, in Proposition \ref{DiffLim} and Definition \ref{soldifference}, taking the diagonal matrix
$v=\frac{1}{u_k I-u_{\hat{k} \hat{k}}}$, and the coefficient matrix
$v\tilde{B}(z)=L_k(z)$, leads to the analytic matrix functions 
\begin{align}\label{Lim:Lp}
    \mathscr{L}^+_k(1+z)^{-1} & \assign
    \left(u_k I-u_{\hat{k} \hat{k}}
    \right)^{1+z} \left( \lim_{p\to\infty}
	p^{-zI-\delta_vB}
	\frac{
    (u_k I-u_{\hat{k} \hat{k}})^{p}
    }{p!}
	\left( \overleftarrow{\prod_{m=1}^{p}}
	L_k(m+z) \right) \right),\quad
    z \notin \mathbb{Z}^{<0}+\Lambda, \\
    \label{Lim:Lm}
    \mathscr{L}^-_k(z) & \assign
    \left( \lim_{p\to\infty}
	\left( \overleftarrow{\prod_{m=-p}^{-1}}
	L_k(m+z)
	\right)
	\frac{
    (u_{\hat{k} \hat{k}}-u_k I)^{p}
    }{p!}
  	p^{zI+\delta_vB} \right) (u_{\hat{k} \hat{k}}-u_k I)^{-z},\quad
    z \notin \mathbb{Z}^{>0}+\Lambda.
\end{align}
and denote the canonical analytic solutions of the system \eqref{Sys:Difference}
as $\mathscr{L}^+_k(z)$, $\mathscr{L}^-_k(z)$, 
with the block matrix indices being the same as those of $L_k(z)$, 
which is $(1,\cdots,k-1,k+1,\cdots,\nu)$.

\begin{thm}\label{Thm:Difeq}
Denote $\arg(u_t-u_s)=-\tau$.
Suppose that none of $u_i$ lies on the segment determined by $u_s$ and $u_t$,
then we have the following results
\begin{subequations}
\begin{align}
    \frac{
	(S_{[\tau]})_{st}
	}{2\pi\mathi} & =
    (u_t-u_s)^{A_{ss}}
	(\mathscr L_t^+(1-A_{tt}^{\bm{r}})^{-1})_{s \hat{t}} \cdot
	\left(\frac{1}{u_tI-u_{\hat{t}\hat{t}}}A_{\hat{t}t}\right),
    \label{StokesDiffeqsol1}\\
    \frac{
	(S_{[\tau]})_{st}
	}{2\pi\mathi}
    & =
    \left(A_{s \hat{s}}\right)\cdot
	\mathscr L_s^-(-A_{ss}^{\bm{l}})_{\hat{s} t}
    (u_t-u_s)^{-A_{tt}}.
    \label{StokesDiffeqsol2}
\end{align}
If we additionally require that none of $u_i$ lies on the line determined by $u_s$ and $u_t$, then we have the following results
\begin{align}
    \label{ConinvDiffeqsol}
    ({C}_{\tau\pm\varepsilon}^{-1})_{\ast s} & =
    \left(I_{\ast \hat{t}}(u_tI_{\hat{t}\hat{t}}-u_{\hat{t}\hat{t}})\right) \cdot
	\mathscr L_t^+(1-A^{\bm{l}})_{\hat{t}s}
    (u_t-u_s)^{-A_{ss}},\\
    \label{ConDiffeqsol}
    ({C}_{\tau\pm\varepsilon})_{t \ast} & =
    (u_t-u_s)^{A_{tt}}
    (\mathscr L_s^-(-A^{\bm{r}})^{-1})_{t\hat{s}} \cdot
	\left(I_{\hat{s}\ast}\right).
\end{align}
\end{subequations}
\end{thm}
\begin{prf}
First, assume that $A_{kk}$ and $A$ are diagonalizable and are all non-resonant.
According to Proposition \ref{DiffLim}, 
the right side of (\ref{StokesDiffeqsol1}-\ref{ConDiffeqsol})
can be interpreted as matrix infinite products \eqref{Lim:Lp} and \eqref{Lim:Lm}
under suitable parameters.
Applying Propositions \ref{FormSolirr} and \ref{FormSolreg}
to Corollary \ref{Cor:LimCoef}, 
and applying \eqref{Lim:Lp} and \eqref{Lim:Lm}
to $v\tilde{B}(z)=L_k(z)$ 
or $v\tilde{B}(z)=z^2L_k(z)^{-1}$,
we obtain \eqref{StokesDiffeqsol1}
and \eqref{ConDiffeqsol}. 
The \eqref{StokesDiffeqsol2}
and \eqref{ConinvDiffeqsol} can be obtained by considering
the conjugate system with coefficient
$-u-\frac{A^\top}{z}$.

Since the dependence of $S_{[\tau]}$ and $C_{\tau\pm\varepsilon}$ on $A$ is analytic and the infinite product on the right hand side of
(\ref{StokesDiffeqsol1}-\ref{ConDiffeqsol}) locally uniformly converges with respect to $A$,  
 the identities
(\ref{StokesDiffeqsol1}-\ref{ConDiffeqsol}) also hold for non-diagonalizable $A$.
Thus, the proof is completed.
\end{prf}

Here the $s$-th row and column 
$(\mathscr{L}_k^+(z)^{-1})_{s \ast}$, 
$(\mathscr{L}_k^-(z))_{\ast s}$
are first defined for
\[
v(t_0)=\diag(v_1I_{n_1},\ldots,t_0v_s I_{n_s},\ldots, v_\mu I_{n_\mu})
\]
with a sufficiently small real number $t_0$ (by Proposition \ref{DiffLim} they are analytic), and then are defined for the original given $v$ through the analytic continuation along $v(t)=\diag(v_1I_{n_1},\ldots,tv_s I_{n_s},\ldots, v_\mu I_{n_\mu})$ by varying $t$ from $t_0$ to 1.
Notice that in this process,
the Stokes matrix on the left side of (\ref{StokesDiffeqsol1}-\ref{StokesDiffeqsol2}) in Theorem \ref{Thm:Difeq} is always analytic, 
ensuring the validity of this process.

\subsection{Bilateral Infinite Matrix Product and Connection Matrices of Difference Systems}

The connection matrix $\mathfrak L_k(z)$ between the canonical analytic solutions $\mathscr{L}^+_k(z)$ and $\mathscr{L}^-_k(z)$
of the differential equation \eqref{kdifference}
can be interpreted as a bilateral infinite matrix product.
The sum of its residues (within one period) is given by the corresponding Stokes matrices. That is
\begin{pro}\label{biinf}
The connection matrix 
$\mathfrak L_k(z) \assign 
 \mathscr{L}^+_k(z)^{-1}
 \mathscr{L}^-_k(z)$
has period $1$,
and the sum of the residues in a period is
\begin{align}
    \tmop{Res} \mathfrak L_k(z)_{s t} =
    (u_k-u_s)^{-A_{ss}}
    \frac{(S_{[\tau]})_{sk}
    (S_{[\tau']})_{kt}}
    {4\pi^2}
    (u_t-u_k)^{A_{tt}},\quad
    k\neq s,t,
\end{align}
where we take the arguments
$\arg(u_k-u_s) = -\tau$,
$\arg(u_t-u_k) = -\tau'$.
\end{pro}

\begin{prf}
Note that we have
\begin{align*}
    \mathfrak L_k(z) = 
    \mathscr{L}^+_k(1+z)^{-1}\cdot
    L_k(z)\cdot
    \mathscr{L}^-_k(z).
\end{align*}
The singularities of $L_k(z)$
are precisely the eigenvalues
$\{-\lambda_1,\ldots,-\lambda_{n_k}\}$
of $-A_{kk}$, 
which provide all the singularities of $\mathfrak L_k(z)$
in a period. 
The non-resonant condition ensures that 
$\mathscr{L}^+_k(1+z)^{-1}$ and 
$\mathscr{L}^-_k(z)$
are analytic at these eigenvalues
$z=-\lambda_i$. 

Without loss of generality, assume that
$A_{kk}$ is semisimple and has a spectral decomposition.
\begin{align*}
    A_{kk} = 
    \lambda_1 P_1 + \cdots +
    \lambda_{n_k} P_{n_k},
\end{align*}
from Theorem \ref{Thm:Difeq} we have
\begin{align*}
    \tmop{Res}_{z=-\lambda_i}
    \mathfrak L_k(z)_{st} & =
    (\mathscr{L}^+_k(1-\lambda_i)^{-1})_{s\hat{k}}\cdot
    \left( -\frac{1}{u_k I - u_{\hat{k} \hat{k}}} 
    A_{\hat{k} k} P_i A_{k \hat{k}}
    \right)\cdot
    \mathscr{L}^-_k(-\lambda_i)_{\hat{k}t} \\
    & = - \left(
    (\mathscr L_k^+(1-A_{kk}^{\bm{r}})^{-1})_{s\hat{k}} \cdot
	\left( \frac{1}{u_kI-u_{\hat{k}\hat{k}}}A_{\hat{k}k}
    \right) \right)P_i \cdot P_i
    \left( \left(A_{k \hat{k}}\right)\cdot
	\mathscr L_k^-(-A_{kk}^{\bm{l}})_{\hat{k} t}
    \right) \\
    & = - (u_k-u_s)^{-A_{ss}}
    \frac{(S_{[\tau]})_{sk}}
    {2\pi\mathi} P_i
    \frac{(S_{[\tau']})_{kt}}
    {2\pi\mathi}
    (u_t-u_k)^{A_{tt}},
\end{align*}
and we finish the proof.
\end{prf}


\section{Path Algebroid}\label{Sect:PathAlg}

We need to understand how to start from $S_d$
to recover $(S_{[\tau]})_{st}$ in Corollary \ref{Cor:LimCoef}
for any fixed $d\notin \tmop{aS}(u)$.
To achieve this, it is necessary to establish the connection between these Stokes matrices. We have already seen in Proposition \ref{Pro:gammajk} how the path class with arguments provides these monodromy data, which naturally leads to the concept of the path algebroid $\mathfrak{S}$ and its representation.
The similar algebroid structure
has also appeared in resurgence theory \cite{mitschi2016divergent}.
Another work of this section is to provide a new interpretation of the action of the braid group on Stokes matrices using the language of the path algebroid.

\subsection{Path Algebroid and its Representation}

\begin{figure}[h!]
    \label{fig:StokesG}
    \centering
    \includegraphics[scale=1]{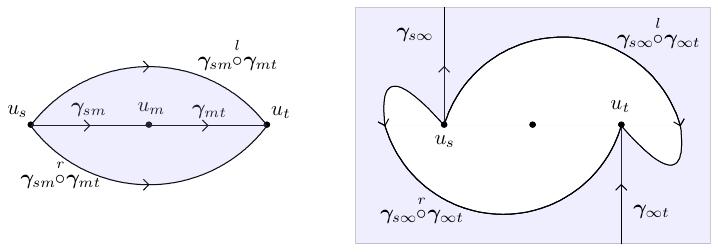}
    \caption{The composition \eqref{Composition} for $m\neq\infty$ and $m=\infty$
    \label{Fig:composition}}
\end{figure}

\begin{defi}[Path algebroid]
The \textbf{path algebroid} $\mathfrak{S}(u_1,\ldots,u_\nu,\infty)$
is an unital algebroid with 
\begin{itemize}
    \item $\nu+1$ objects $u_1,\ldots,u_\nu\in\mathbb{C},u_\infty\assign\infty$;
    
    \item For any two objects $u_s,u_t$, the hom-set $\operatorname{Hom}(u_s,u_t)$
    is the complex linear space generated by 
    $\delta_{s t} \assign \left\{\begin{array}{ll}
     1_s, & s = t\\
     0, & s \neq t
    \end{array}\right.$
    and all path classes with arguments in  
    $\mathbb{C}\setminus\{u_1,\ldots,u_\nu\}$
    from $u_s$ to $u_t$.
    For any morphisms
    $\boldsymbol{\gamma}_{st}\in \operatorname{Hom}(u_s,u_t)$,
    we have
    \begin{align*}
        1_s \boldsymbol{\gamma}_{st} = \boldsymbol{\gamma}_{st} =
        \boldsymbol{\gamma}_{st} 1_t.
    \end{align*}
    When there is no ambiguity, we will denote $1_s$ as $1$.
    
    \item For any two morphisms
    $\boldsymbol{\gamma}_{sm}\in \operatorname{Hom}(u_s,u_m)$,
    $\boldsymbol{\gamma}_{mt}\in \operatorname{Hom}(u_m,u_t)$,
    the multiplication is defined by
    \begin{align}\label{Eq:StokesW}
    \boldsymbol{\gamma}_{sm}
    \boldsymbol{\gamma}_{mt} \assign
    \boldsymbol{\gamma}_{sm} \overset{r}{\circ}
    \boldsymbol{\gamma}_{mt} -
    \boldsymbol{\gamma}_{sm} \overset{l}{\circ}
    \boldsymbol{\gamma}_{mt}.
    \end{align}
    For any path class with arguments $\boldsymbol{\gamma}_{sm}$
    from $u_s$ to $u_m$,
    and  $\boldsymbol{\gamma}_{mt}$
    from $u_m$ to $u_t$,
    take the representative elements 
    $\tilde{\gamma}_{sm}=(\gamma_{sm},\theta_s,\theta_m)$ and 
    $\tilde{\gamma}_{mt}=(\gamma_{mt},\theta_m,\theta_t)$.
    Denote $\gamma_{sm}\overset{l}{\circ}\gamma_{mt},\gamma_{sm}\overset{r}{\circ}\gamma_{mt}$ 
    as paths in $\mathbb{C}\setminus\{u_1,\ldots,u_\nu,\infty\}$,
    each obtained by deviating $\gamma_{sm}\circ\gamma_{mt}$ to the left and right at $u_m$, respectively (see Figure \ref{Fig:composition}).
    \begin{alignat}{2}
    \tilde{\gamma}_{sm} \overset{l}{\circ}
    \tilde{\gamma}_{mt} & \assign
    (\gamma_{sm}\overset{l}{\circ}\gamma_{mt},\theta_s,\theta_t),&
    \quad
    \tilde{\gamma}_{sm} \overset{r}{\circ}
    \tilde{\gamma}_{mt} & \assign
    (\gamma_{sm}\overset{r}{\circ}\gamma_{mt},\theta_s,\theta_t),\\
    \label{Composition}
    \boldsymbol{\gamma}_{sm} \overset{l}{\circ}
    \boldsymbol{\gamma}_{mt} & \assign
    (\tilde{\gamma}_{sm} \overset{l}{\circ}
    \tilde{\gamma}_{mt})/\sim,&
    \quad
    \boldsymbol{\gamma}_{sm} \overset{r}{\circ}
    \boldsymbol{\gamma}_{mt} & \assign
    (\tilde{\gamma}_{sm} \overset{r}{\circ}
    \tilde{\gamma}_{mt})/\sim.
    \end{alignat}
\end{itemize}
\end{defi}
It can be seen that
the composition $\overset{l}{\circ}, \overset{r}{\circ}$
is well-defined, associative and compatible with each other, i.e.
    \begin{align}
        (\boldsymbol{\gamma}_{sm_1} \overset{l}{\circ}
        \boldsymbol{\gamma}_{m_1m_2})\overset{r}{\circ}
        \boldsymbol{\gamma}_{m_2t} =
        \boldsymbol{\gamma}_{sm_1} \overset{l}{\circ}
        (\boldsymbol{\gamma}_{m_1m_2} \overset{r}{\circ}
        \boldsymbol{\gamma}_{m_2t}),\quad
        (\boldsymbol{\gamma}_{sm_1} \overset{r}{\circ}
        \boldsymbol{\gamma}_{m_1m_2}) \overset{l}{\circ}
        \boldsymbol{\gamma}_{m_2t} =
        \boldsymbol{\gamma}_{sm_1} \overset{r}{\circ}
        (\boldsymbol{\gamma}_{m_1m_2} \overset{l}{\circ}
        \boldsymbol{\gamma}_{m_2t}),
    \end{align}
    Therefore, the multiplication defined in \eqref{Eq:StokesW} is also associative.
    Moreover, the composition $\overset{l}{\circ}, \overset{r}{\circ}$ 
    each has $\boldsymbol{e}_s^{l}$ and $\boldsymbol{e}_s^{r}$ 
    as its identity elements, 
    and any path class with arguments from $u_s$ to $u_t(t\neq\infty)$
    has the following inverse
    \begin{align}
        (\gamma_{st},\theta_s,\theta_t)^{\overset{l}{\circ}-1} =
        (\gamma_{st}^{\tmop{op}},\theta_t+\pi,\theta_s-\pi),\quad
        (\gamma_{st},\theta_s,\theta_t)^{\overset{r}{\circ}-1} =
        (\gamma_{st}^{\tmop{op}},\theta_t-\pi,\theta_s+\pi).
    \end{align}
    Therefore, these path classes with arguments even form a groupoid.
\begin{defi}
A representation $(\mathcal{S},V_1,\ldots,V_\nu)$ of the path algebroid $\mathfrak{S}(u_1,\ldots,u_\nu,\infty)$ refers to a family of complex vector spaces $\{V_s:s=1,\ldots,\nu,\infty\}$ with $V_\infty=V_1\oplus\cdots\oplus V_\nu$, and morphism $\mathcal{S}$ satisfying
\begin{align}
    \mathcal{S}(u_s) = V_s,\quad
    \operatorname{Hom}(u_s,u_t)
    \xrightarrow{\mathcal{S}}
    \operatorname{Hom}(V_s,V_t),\quad
    \mathcal{S}(\boldsymbol{\gamma}_{sm}
    \boldsymbol{\gamma}_{mt}) =
    \mathcal{S}(\boldsymbol{\gamma}_{sm})
    \mathcal{S}(\boldsymbol{\gamma}_{mt}).
\end{align}
\end{defi}


\begin{rmk}\label{Rmk:SW}
    Note that we have
    \begin{subequations}
    \begin{align}\label{erCompL}
        (\gamma_{st},\theta_s,\theta_t) \overset{l}{\circ} 
        \boldsymbol{e}_t^{r} & =
        (\gamma_{st},\theta_s,\theta_t+2\pi),\quad
        \boldsymbol{e}_s^{r} \overset{l}{\circ} 
        (\gamma_{st},\theta_s,\theta_t) =
        (\gamma_{st},\theta_s-2\pi,\theta_t),\\
        \label{elCompR}
        (\gamma_{st},\theta_s,\theta_t) \overset{r}{\circ} 
        \boldsymbol{e}_t^{l} & =
        (\gamma_{st},\theta_s,\theta_t-2\pi),\quad
        \boldsymbol{e}_s^{l} \overset{r}{\circ} 
        (\gamma_{st},\theta_s,\theta_t) =
        (\gamma_{st},\theta_s+2\pi,\theta_t),
    \end{align}
    \end{subequations}
    provided that when the set $G$
    contains at least one path class with arguments 
    between any two distinct objects $u_s$ and $u_t$, 
    then $G$ can generate the entire path algebroid 
    together with
    $\{\boldsymbol{e}_s^{l},\boldsymbol{e}_s^{r}\}_{s=1,\ldots,\nu,\infty}$
    and an additional relation
    $(1+\boldsymbol{e}_s^{l})^{-1} =
    1-\boldsymbol{e}_s^{r}$,
    under the multiplication \eqref{Eq:StokesW}.
\end{rmk}

\begin{thm}\label{Thm:SWeight}
$(\mathcal{S}_{u,A},V_1,\ldots,V_\nu)$ is a representation of the path algebroid
$\mathfrak{S}(u_1,\ldots,u_\nu,\infty)$,
where $V_i$ is the eigenspace corresponding to the eigenvalue $u_i$ of $u$.
\end{thm}

\begin{prf}
We have the set of generators 
$G=\{
\boldsymbol{e}_s^{l},\boldsymbol{e}_s^{r},
\boldsymbol{\gamma}_{st}^{(d)}:
s,t=1,\ldots,\nu,\infty;s\neq t\}$. 
It suffices to prove that for any path class with arguments $\boldsymbol{\gamma}_{sm}$
and generator $\boldsymbol{\gamma}_{mt}\in G$,
we have $\mathcal{S}(\boldsymbol{\gamma}_{sm}\boldsymbol{\gamma}_{mt}) = 
\mathcal{S}(\boldsymbol{\gamma}_{sm})\mathcal{S}(\boldsymbol{\gamma}_{mt})$.

For $m=t$,
according to 
\eqref{erCompL}, \eqref{elCompR} and 
\eqref{LoopWeight} in Proposition \ref{Pro:gammajk},
we can first verify that
$\mathcal{S}(1-\boldsymbol{e}_s^{r}) =
 \mathe^{-2\pi\mathi A_{ss}} =
 \mathcal{S}(1+\boldsymbol{e}_s^{l})^{-1}$ 
and according to Definition \ref{Def:WeightS} we have
\begin{align*}
    \mathcal{S}(\boldsymbol{\gamma}_{sm}\boldsymbol{e}_m^{l}) & =
    \mathcal{S}(\boldsymbol{\gamma}_{sm})
    (\mathe^{2\pi\mathi A_{mm}} - I) =
    \mathcal{S}(\boldsymbol{\gamma}_{sm})
    \mathcal{S}(\boldsymbol{e}_m^{l}),\\
    \mathcal{S}(\boldsymbol{\gamma}_{sm}\boldsymbol{e}_m^{r}) & =
    \mathcal{S}(\boldsymbol{\gamma}_{sm})
    (I - \mathe^{-2\pi\mathi A_{mm}}) =
    \mathcal{S}(\boldsymbol{\gamma}_{sm})
    \mathcal{S}(\boldsymbol{e}_m^{r}),
\end{align*}
where $A_{\infty\infty}\assign A$.

Next, we assume that $m \neq t$.
For $t\neq \infty$,
according to \eqref{Weight:S}, \eqref{Weight:Sinf}
and Corollary \ref{ProveWeight} 
we have
\begin{align*}
    \mathcal{S}(\boldsymbol{\gamma}_{sm} \overset{r}{\circ}
    \boldsymbol{\gamma}_{mt}^{(d)}) -
    \mathcal{S}(\boldsymbol{\gamma}_{sm} \overset{l}{\circ}
    \boldsymbol{\gamma}_{mt}^{(d)}) 
    & =
    \left\{
    \begin{array}{ll}
    \mathcal{S}(\boldsymbol{\gamma}_{sm})
    (S_d^+ -
    \mathe^{- 2\pi\mathi\delta_u A} S_d^-)_{mt}, 
    & m\neq\infty,
      \tmop{Im} (u_{m} \mathe^{\mathi d}) >
      \tmop{Im} (u_{t} \mathe^{\mathi d})\\
    \mathcal{S}(\boldsymbol{\gamma}_{sm})
    (\mathe^{2\pi\mathi\delta_u A}S_d^+ -
     S_d^-)_{mt},
    & m\neq\infty,
      \tmop{Im} (u_{t} \mathe^{\mathi d}) >
      \tmop{Im} (u_{m} \mathe^{\mathi d})\\
    \mathcal{S}(\boldsymbol{\gamma}_{sm})
    ({\rm e}^{\pi{\rm i} A} C_d^{-1} (I-U_d^{[\infty]})
    )_{\ast t},
    & m=\infty
    \end{array}
    \right.\\
    & =
    \mathcal{S}(\boldsymbol{\gamma}_{sm})
    \mathcal{S}(\boldsymbol{\gamma}_{mt}^{(d)}).
\end{align*}
For $t=\infty$,
according to \eqref{Weight:O}, \eqref{Weight:H0p}
and Corollary \ref{ProveWeight} 
we have
\begin{align*}
    \mathcal{S}(\boldsymbol{\gamma}_{sm} \overset{r}{\circ}
    \boldsymbol{\gamma}_{mt}^{(d)}) -
    \mathcal{S}(\boldsymbol{\gamma}_{sm} \overset{l}{\circ}
    \boldsymbol{\gamma}_{mt}^{(d)}) 
    & =
    \mathcal{S}(\boldsymbol{\gamma}_{sm})
    ((I-U_{d-\pi}^{[m]}) C_{d-\pi} {\rm e}^{\pi{\rm i}A})_{m\ast}\\
    & =
    \mathcal{S}(\boldsymbol{\gamma}_{sm})
    (E_{mm}S_{d-\pi}^+(I-{\rm e}^{-2\pi{\rm i}M_d})
    C_{d-\pi} {\rm e}^{\pi{\rm i}A})_{m\ast}\\
    & =
    \mathcal{S}(\boldsymbol{\gamma}_{sm})
    \mathcal{S}(\boldsymbol{\gamma}_{mt}^{(d)}).
\end{align*}
\noindent
Thus, we complete the proof.
\end{prf}\\
\noindent
In particular, starting from $\mathe^{2\pi\mathi\delta_uA}$
    and $S_d$ in a specific direction $d$, 
    or $(S_{[\tau]})_{st}$ in Corollary \ref{Cor:LimCoef}, 
    is sufficient to recover Stokes matrices in all directions.



\subsection{Stokes Matrices and Quantum Groups}
Let us first recall the quantum group relations of the Stokes matrices. 

Let us take the standard $R$-matrix
$R\in {\rm End}(\mathbb{C}^\nu)\otimes {\rm End}(\mathbb{C}^\nu)$ with $q=\mathe^{\pi\mathi\hbar}\in\mathbb{C}$, see e.g., \cite{Jimbo2}\cite{Faddeev1988},
\begin{equation}\label{sRmatrix}
R=\sum_{i\ne j, i,j=1}^n E_{ii}\otimes E_{jj}+q\sum_{i=1}^n E_{ii}\otimes E_{ii}+(q-q^{-1})\sum_{1\le j<i\le n}E_{ij}\otimes E_{ji}.
\end{equation} 

\begin{thm}\cite{Xu2}\label{qStokes}
For any $h\notin\mathbb{Q}$ and distinct $u_1,...,u_\nu$, 
consider the quantum confluent hypergeometric equation \eqref{introqeq} associated to a representation $V$. Assume that the direction $d\notin \tmop{aS}(u)$ is such that
    \begin{equation*}
      \tmop{Im} (u_{1} \mathe^{\mathi d}) >
      \tmop{Im} (u_{2} \mathe^{\mathi d}) > \cdots > 
      \tmop{Im} (u_{\nu}  \mathe^{\mathi d}),
    \end{equation*} then the (modified) Stokes matrices 
$L_{\pm} = q^{\mp \delta_u \mathbf{E}} 
S_d^{\pm} (u, - \hbar \mathbf{E})\in {\rm End}(V)\otimes{\rm End}(\mathbb{C}^\nu)$, of the quantum confluent hypergeometric equation \eqref{introqeq} associated to a representation $V$, satisfy
\begin{equation}\label{RLL=LLR}
R^{12}L_\pm^{(1)}L_\pm^{(2)}=L_\pm^{(2)}L_\pm^{(1)} R^{12},
\end{equation}
\begin{equation}\label{RLpLm=LmLpR}
R^{12}L_+^{(1)}L_-^{(2)} =L^{(2)}_{+}L^{(1)}_{-}(u) R^{12}.
\end{equation}
Here the convention is that if we write
$L_{\pm}=\sum_{i,j}(L_\pm)_{ij}\otimes E_{ij}\in {\rm End}(V)\otimes {\rm End}(\mathbb{C}^\nu)$,
then 
\[
L_{\pm}^{(1)}:= \sum_{i,j}(L_\pm)_{ij}\otimes E_{ij}\otimes I,\quad
L_\pm^{(2)}:= \sum_{i,j}(L_\pm)_{ij}\otimes I\otimes E_{ij},\quad
\text{and } R^{12}:=I\otimes R
\] 
as elements in ${\rm End}(V)\otimes  {\rm End}(\mathbb{C}^\nu)\otimes {\rm End}(\mathbb{C}^\nu)$.
\end{thm}

Denote $\mathcal{s}_{i j}^\pm = S_d^{\pm} (u, - \hbar \mathbf{E})_{ij}$ and $\mathcal{h}_i=e_{ii}$.
By directly comparing both sides of 
\eqref{RLL=LLR}, \eqref{RLpLm=LmLpR},
we have the commutation relations between
$q^{\mathcal{h}_k}$ and other elements:
\begin{subequations}
\begin{align}\label{ComRelBegin}
  q^{\mathcal{h}_i} \cdot q^{\mathcal{h}_j} & = q^{\mathcal{h}_j} \cdot
  q^{\mathcal{h}_i}, \\
  q^{\mathcal{h}_k} \cdot \mathcal{s}^{\pm}_{i j} & = \mathcal{s}^{\pm}_{ij} \cdot 
  q^{\mathcal{h}_k + \delta_{i k} - \delta_{k j}},
\end{align}   
\end{subequations}
as well as the commutation relations when the indices $i_1,j_1,i_2,j_2$ are distinct:
\begin{subequations}
\begin{align}
\frac{
\mathcal{s}^+_{i_1 j_1} 
\mathcal{s}^+_{i_2 j_2} -
\mathcal{s}^+_{i_2 j_2} 
\mathcal{s}^+_{i_1 j_1}}{q - q^{- 1}} 
& = \left\{\begin{array}{ll}
    \mathcal{s}^+_{i_2 j_1} \cdot 
    \mathcal{s}^+_{i_1 j_2}, & i_2 < i_1 < j_2 < j_1\\
    -\mathcal{s}^+_{i_2 j_1} \cdot 
    \mathcal{s}^+_{i_1 j_2}, & i_1 < i_2 < j_1 < j_2\\
    \multirow{4}{*}{$0,$}      & i_2 < i_1 < j_1 < j_2\\
                             & i_1 < i_2 < j_2 < j_1\\
                             & i_1 < j_1 < i_2 < j_2\\
                             & i_2 < j_2 < i_1 < j_1
    \end{array}\right., \\
\frac{
\mathcal{s}^-_{i_1 j_1} 
\mathcal{s}^-_{i_2 j_2} -
\mathcal{s}^-_{i_2 j_2} 
\mathcal{s}^-_{i_1 j_1}}{q - q^{- 1}} 
& = \left\{\begin{array}{ll}
    \mathcal{s}^-_{i_2 j_1} \cdot 
    \mathcal{s}^-_{i_1 j_2}, & j_2 < j_1 < i_2 < i_1\\
    -\mathcal{s}^-_{i_2 j_1} \cdot 
    \mathcal{s}^-_{i_1 j_2}, & j_1 < j_2 < i_1 < i_2\\
    \multirow{4}{*}{$0,$}      & j_1 < j_2 < i_2 < i_1\\
                             & j_2 < j_1 < i_1 < i_2\\
                             & j_2 < i_2 < j_1 < i_1\\
                             & j_1 < i_1 < j_2 < i_2
    \end{array}\right., \\
\frac{
\mathcal{s}^+_{i_1 j_1} 
\mathcal{s}^-_{i_2 j_2} -
\mathcal{s}^-_{i_2 j_2} 
\mathcal{s}^+_{i_1 j_1}}{q - q^{- 1}} 
& = \left\{\begin{array}{ll}
    \mathcal{s}^-_{i_2 j_1} \cdot 
    \mathcal{s}^+_{i_1 j_2}, & i_1 < j_2 < j_1 < i_2\\
    - q^{2 (\mathcal{h}_{i_1} -\mathcal{h}_{i_2})} 
    \mathcal{s}^+_{i_2 j_1} \cdot
    \mathcal{s}^-_{i_1 j_2}, & j_2 < i_1 < i_2 < j_1\\
    \multirow{4}{*}{$0,$}      & i_1 < j_1 < j_2 < i_2\\
                             & j_2 < i_2 < i_1 < j_1\\
                             & i_1 < j_2 < i_2 < j_1\\
                             & j_2 < i_1 < j_1 < i_2
    \end{array}\right.,
\end{align}
\end{subequations}
and the commutation relations when some of the indices coincide:
\begin{subequations}
\begin{align}
\frac{ q \cdot 
\mathcal{s}^+_{i_1 j_1} 
\mathcal{s}^+_{i_2 j_2} -
\mathcal{s}^+_{i_2 j_2} 
\mathcal{s}^+_{i_1 j_1}}{q - q^{- 1}} 
& = \left\{\begin{array}{ll}
    q \cdot 
    \mathcal{s}^+_{i_2 j_1}, & i_2 < i_1 = j_2 < j_1\\
    \mathcal{s}^+_{i_2 j_1} \cdot 
    \mathcal{s}^+_{i_1 j_2}, & i_2 = i_1 < j_2 < j_1\\
    \mathcal{s}^+_{i_1 j_2} \cdot 
    \mathcal{s}^+_{i_2 j_1}, & i_2 < i_1 < j_2 = j_1\\
    \multirow{2}{*}{$0,$}      & i_1 = i_2 < j_1 < j_2\\
                             & i_1 < i_2 < j_1 = j_2
    \end{array}\right., \\
\frac{ q \cdot 
\mathcal{s}^-_{i_1 j_1} 
\mathcal{s}^-_{i_2 j_2} -
\mathcal{s}^-_{i_2 j_2} 
\mathcal{s}^-_{i_1 j_1}}{q - q^{- 1}} 
& = \left\{\begin{array}{ll}
    \mathcal{s}^-_{i_1 j_2}, & j_2 < j_1 = i_2 < i_1\\
    \mathcal{s}^-_{i_2 j_1} \cdot 
    \mathcal{s}^-_{i_1 j_2}, & j_2 < j_1 < i_2 = i_1\\
    \mathcal{s}^-_{i_1 j_2} \cdot 
    \mathcal{s}^-_{i_2 j_1}, & j_2 = j_1 < i_2 < i_1\\
    \multirow{2}{*}{$0,$}      & j_1 < j_2 < i_1 = i_2\\
                             & j_1 = j_2 < i_1 < i_2
    \end{array}\right., \\
\frac{ q \cdot 
\mathcal{s}^+_{i_1 j_1} 
\mathcal{s}^-_{i_2 j_2} -
\mathcal{s}^-_{i_2 j_2} 
\mathcal{s}^+_{i_1 j_1}}{q - q^{- 1}} 
& = \left\{\begin{array}{ll}
    q \cdot 
    \mathcal{s}^-_{i_2 j_1}, & i_1 = j_2 < j_1 < i_2\\
    \mathcal{s}^+_{i_1 j_2}, & i_1 < j_2 < j_1 = i_2\\
    \multirow{2}{*}{$
    \mathcal{s}^+_{i_1 j_1} \cdot 
    \mathcal{s}^-_{i_2 j_2},$}& i_1 < j_2 = j_1 < i_2\\
                             & j_2 < i_1 = i_2 < j_1\\
    - q \cdot q^{2 (\mathcal{h}_{i_1} -\mathcal{h}_{i_2})} 
    \mathcal{s}^+_{i_2 j_1}, & j_2 = i_1 < i_2 < j_1\\
    - q^{2 (\mathcal{h}_{i_1} -\mathcal{h}_{i_2})} 
    \mathcal{s}^-_{i_1 j_2}, & j_2 < i_1 < i_2 = j_1
  \end{array}\right., \\
\frac{ q \cdot 
\mathcal{s}^+_{i j} 
\mathcal{s}^-_{j i} - 
q^{- 1} \cdot
\mathcal{s}^-_{j i} 
\mathcal{s}^+_{i j}}{q - q^{- 1}} 
& = 1 - q^{2(\mathcal{h}_i -\mathcal{h}_j)}, \quad i < j. 
\label{ComRelEnd}\end{align}
\end{subequations}
Equivalently, the above theorem can be stated as follows.
\begin{thm}\label{Squantumgroup}\cite{Xu2}
The Stokes matrices of the equation \ref{introqeq}
induce a representation of $U_q(\mathfrak{gl}_\nu)$ on $V$ via the map
\begin{align}
    q^{\pm h_j}\mapsto q^{\pm\mathcal{h}_j},\quad
    f_i\mapsto \frac{1}{q-q^{-1}} \mathcal{s}^+_{i,i+1},\quad
    e_i\mapsto-\frac{1}{q-q^{-1}} q^{h_{i+1}}\mathcal{s}^-_{i+1,i}q^{-h_i}
    \in {\rm End}(V).
\end{align}
Here recall that the Drinfeld-Jimbo quantum group $U_q(\mathfrak{gl}_\nu)$ is an associative algebra with
generators $\{f_i,e_i,q^{\pm h_j}:1\leqslant i\leqslant \nu-1, 1\leqslant j \leqslant \nu\}$ and relations 
\cite{ding1993isomorphism}
\begin{gather}
    q^{-h_j}=(q^{h_j})^{-1},\quad
    q^{h_j}f_iq^{-h_j}=q^{\delta_{ij}-\delta_{j,i+1}}f_i,\quad
    q^{h_j}e_iq^{-h_j}=q^{\delta_{i+1,j}-\delta_{ji}}e_i,\\
    f_ie_i-e_if_i=\frac{q^{h_{i}-h_{i+1}}-q^{h_{i+1}-h_{i}}}{q-q^{-1}},\quad
    f_{i_2}e_{i_1}-e_{i_1}f_{i_2}=0,\quad 
    i_1\neq i_2,\\
    f_{i_1}^2f_{i_2}-(q+q^{-1})f_{i_1}f_{i_2}f_{i_1}+f_{i_2}f_{i_1}^2=0,\quad
    |i_1-i_2|=1;\quad
    f_{i_1}f_{i_2} = f_{i_2}f_{i_1},\quad
    |i_1-i_2|\geqslant 2,\\
    e_{i_1}^2e_{i_2}-(q+q^{-1})e_{i_1}e_{i_2}e_{i_1}+e_{i_2}e_{i_1}^2=0,\quad
    |i_1-i_2|=1;\quad
    e_{i_1}e_{i_2} = e_{i_2}e_{i_1},\quad
    |i_1-i_2|\geqslant 2.
\end{gather}
\end{thm}

\subsection{A Morphism from the Algebroid \texorpdfstring{$\mathfrak{S}(u_1,\ldots,u_\nu)$}{S} to the Quantum Group \texorpdfstring{$U_q(\mathfrak{gl}_\nu)$}{U}}

Theorem \ref{Thm:SWeight} and Theorem \ref{qStokes}
indicate that there exists a morphism $\pi$ from
$\mathfrak{S}(u_1,\ldots,u_\nu)$ 
to $U_q(\mathfrak{gl}_\nu)$.
\begin{cor}
    Suppose that for direction $d$ we have
    \begin{equation*}
      \tmop{Im} (u_{1} \mathe^{\mathi d}) >
      \tmop{Im} (u_{2} \mathe^{\mathi d}) > \cdots > 
      \tmop{Im} (u_{\nu}  \mathe^{\mathi d}),
    \end{equation*}
    then there is a morphsim 
    $\pi:\mathfrak{S}(u_1,\ldots,u_\nu)\to U_q(\mathfrak{gl}_\nu)$,
    induced by
    \begin{align*}
        \pi(\boldsymbol{\gamma}_{i,i+1}^{(d)}) = 
        (q-q^{-1})f_i,\quad
        \pi(\boldsymbol{\gamma}_{i+1,i}^{(d)}) = 
        (q-q^{-1})q^{-h_{i+1}}e_iq^{h_i},\quad
        \pi(\boldsymbol{e}_{i}^{l}) = q^{-2h_i}-1,\quad
        \pi(\boldsymbol{e}_{i}^{r}) = 1-q^{2h_i}.
    \end{align*}
\end{cor}
Moreover, equations (\ref{ComRelBegin}-\ref{ComRelEnd})
provide the commutation relations for each pair of generators.
It is natural to consider the explicit commutation relation
for other pairs of path classes with arguments under this morphsim. In order to to represent the commutation relation
for the representative elements
$\tilde{\gamma}_{s_1t_1} = 
(\gamma_{s_1t_1},\theta_s^{(1)},\theta_t^{(1)})$,
$\tilde{\gamma}_{s_2t_2} = 
(\gamma_{s_2t_2},\theta_s^{(2)},\theta_t^{(2)})$ 
of path classes with arguments,
we need to introduce the following notation 
\begin{align}\nonumber
    \tilde{\gamma}_{s_1t_1} \wedge
    \tilde{\gamma}_{s_2t_2} & \assign
    \sum_{m\in \gamma_{s_1t_1}\cap \gamma_{s_2t_2} }
    \operatorname{sgn}_{m}
    (\gamma_{s_1t_1},\gamma_{s_2t_2})
    \left( \gamma_{s_1t_1} \overset{m}{\ast} \gamma_{s_2t_2},
    \theta_s^{(1)}, \theta_t^{(2)} \right) 
    \left( \gamma_{s_2t_2} \overset{m}{\ast} \gamma_{s_1t_1},
    \theta_s^{(2)}, \theta_t^{(1)} \right)\\
    & \phantom{=}\quad
    - \delta_{t_1 s_2}
    \begin{array}{l}
        \delta_{ \theta_t^{(1)} > \theta_s^{(2)} } 
        \tilde{\gamma}_{s_1t_1} \overset{r}{\circ} 
        \tilde{\gamma}_{s_2t_2}\\
        \delta_{ \theta_t^{(1)} < \theta_s^{(2)} } 
        \tilde{\gamma}_{s_1t_1} \overset{l}{\circ} 
        \tilde{\gamma}_{s_2t_2}
    \end{array}
    + \delta_{t_2 s_1} 
    \begin{array}{l}
        \delta_{ \theta_t^{(2)} > \theta_s^{(1)} } 
        \tilde{\gamma}_{s_2t_2} \overset{r}{\circ} 
        \tilde{\gamma}_{s_1t_1}\\
        \delta_{ \theta_t^{(2)} < \theta_s^{(1)} }
        \tilde{\gamma}_{s_2t_2} \overset{l}{\circ} 
        \tilde{\gamma}_{s_1t_1}
    \end{array},\label{wedge}
\end{align}
where 
$| \theta_t^{(1)} - \theta_s^{(2)} |,
 | \theta_t^{(2)} - \theta_s^{(1)} | < \pi$,
and the intersection point $m$ does not include endpoints.
Here $\tmop{sgn}_m$ denotes the orientation
of $\gamma_1$ and $\gamma_2$
at the intersection point $m$.
\begin{align}
    \tmop{sgn}_m (\gamma_{s_1t_1},\gamma_{s_2t_2}) & \assign 
    \left\{\begin{array}{ll}
    1, &
    \vcenter{\hbox{\includegraphics[width=0.135\textwidth]{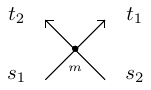}}} \\
    - 1, &
    \vcenter{\hbox{\includegraphics[width=0.135\textwidth]{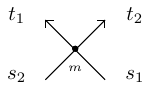}}}
  \end{array}\right.
\end{align}
And $\gamma_{s_1t_1} \overset{m}{\ast} \gamma_{s_2t_2}$
represents a path that starts along path $\gamma_{s_1t_1}$ and, 
after reaching the intersection point $m$,
turns onto path $\gamma_{s_2t_2}$.
For example, 
for distinct $u_s,u_m,u_t$
and $u_m\notin [u_s,u_t]$,
we have
\begin{align*}
    [u_s,u_m]\wedge [u_m,u_t] = - [u_s,u_t].
\end{align*}

\begin{lem}\label{Lem:Leb}
The wedge
defined in \eqref{wedge} satisfies 
the right Leibniz law.
For the representative elements of path class with arguments,
\begin{align*}
    \tilde{\gamma}_{s_1m} =
    (\gamma_{s_1m}, \theta_s^{(1)}, \theta_m),\quad
    \tilde{\gamma}_{mt_1} =
    (\gamma_{mt_1}, \theta_m, \theta_t^{(1)}),\quad
    \tilde{\gamma}_{s_2t_2} =
    (\gamma_{s_2t_2}, \theta_s^{(2)}, \theta_t^{(2)}),
\end{align*}
we have
\begin{align}\nonumber
    (\tilde{\gamma}_{s_1m}
     \tilde{\gamma}_{mt_1})
    \wedge \tilde{\gamma}_{s_2t_2} & \assign
    (\tilde{\gamma}_{s_1m} \overset{r}{\circ} 
    \tilde{\gamma}_{mt_1})
    \wedge \tilde{\gamma}_{s_2t_2} -
    (\tilde{\gamma}_{s_1m} \overset{l}{\circ} 
    \tilde{\gamma}_{mt_1})
    \wedge \tilde{\gamma}_{s_2t_2}\\
    & =
    (\tilde{\gamma}_{s_1m}
     \wedge \tilde{\gamma}_{s_2t_2} )
     \tilde{\gamma}_{mt_1} +
     \tilde{\gamma}_{s_1m}
    (\tilde{\gamma}_{mt_1}
     \wedge \tilde{\gamma}_{s_2t_2} ).
\end{align}
\end{lem}

\begin{prf}
Note that we have
\begin{eqnarray*}
  \left( \tilde{\gamma}_{s_1 m} \overset{r}{\circ} \tilde{\gamma}_{m t_1}
  \right) \wedge \tilde{\gamma}_{s_2 t_2} & = & \sum_{p \in \gamma_{s_1 m}
  \cap \gamma_{s_2 t_2}} \tmop{sgn}_p (\gamma_{s_1 m}, \gamma_{s_2 t_2})
  \left( \tilde{\gamma}_{s_1 m} \overset{p}{\ast} \tilde{\gamma}_{s_2 t_2}
  \right) \left( \left( \tilde{\gamma}_{s_2 t_2} \overset{p}{\ast}
  \tilde{\gamma}_{s_1 m} \right) \overset{r}{\circ} \tilde{\gamma}_{m t_1}
  \right)\\
  &  & + \sum_{p \in \gamma_{m t_1} \cap \gamma_{s_2 t_2}} \tmop{sgn}_p
  (\gamma_{m t_1}, \gamma_{s_2 t_2}) \left( \tilde{\gamma}_{s_1 m}
  \overset{r}{\circ} \left( \tilde{\gamma}_{m t_1} \overset{p}{\ast}
  \tilde{\gamma}_{s_2 t_2} \right) \right) \left( \tilde{\gamma}_{s_2 t_2}
  \overset{p}{\ast} \tilde{\gamma}_{m t_1} \right)\\
  &  & + \left\{\begin{array}{ll}
    \tilde{\gamma}_{s_1 m} \left( \tilde{\gamma}_{s_2 t_2} \overset{r}{\circ}
    \tilde{\gamma}_{m t_1} \right), & t_2 = m, \theta_m < \theta_t^{(2)}\\
    0, & \text{etc.}
  \end{array}\right.\\
  &  & - \left\{\begin{array}{ll}
    \left( \tilde{\gamma}_{s_1 m} \overset{r}{\circ} \tilde{\gamma}_{s_2 t_2}
    \right) \tilde{\gamma}_{m t_1}, & m = s_2, \theta_s^{(2)} < \theta_m\\
    0, & \text{etc.}
  \end{array}\right.\\
  &  & + \delta_{t_2 s_1} \begin{array}{c}
    \delta_{\theta_s^{(1)} < \theta_t^{(2)}} \tilde{\gamma}_{s_2 t_2}
    \overset{r}{\circ} \tilde{\gamma}_{s_1 m} \overset{r}{\circ}
    \tilde{\gamma}_{m t_1}\\
    \delta_{\theta_s^{(1)} > \theta_t^{(2)}} \tilde{\gamma}_{s_2 t_2}
    \overset{l}{\circ} \tilde{\gamma}_{s_1 m} \overset{r}{\circ}
    \tilde{\gamma}_{m t_1}
  \end{array} - \delta_{t_1 s_2} \begin{array}{c}
    \delta_{\theta_s^{(2)} < \theta_t^{(1)}} \tilde{\gamma}_{s_1 m}
    \overset{r}{\circ} \tilde{\gamma}_{m t_1} \overset{r}{\circ}
    \tilde{\gamma}_{s_2 t_2}\\
    \delta_{\theta_s^{(2)} > \theta_t^{(1)}} \tilde{\gamma}_{s_1 m}
    \overset{r}{\circ} \tilde{\gamma}_{m t_1} \overset{l}{\circ}
    \tilde{\gamma}_{s_2 t_2}
  \end{array}\\
  \left( \tilde{\gamma}_{s_1 m} \overset{l}{\circ} \tilde{\gamma}_{m t_1}
  \right) \wedge \tilde{\gamma}_{s_2 t_2} & = & \sum_{p \in \gamma_{s_1 m}
  \cap \gamma_{s_2 t_2}} \tmop{sgn}_p (\gamma_{s_1 m}, \gamma_{s_2 t_2})
  \left( \tilde{\gamma}_{s_1 m} \overset{p}{\ast} \tilde{\gamma}_{s_2 t_2}
  \right) \left( \left( \tilde{\gamma}_{s_2 t_2} \overset{p}{\ast}
  \tilde{\gamma}_{s_1 m} \right) \overset{l}{\circ} \tilde{\gamma}_{m t_1}
  \right)\\
  &  & + \sum_{p \in \gamma_{m t_1} \cap \gamma_{s_2 t_2}} \tmop{sgn}_p
  (\gamma_{m t_1}, \gamma_{s_2 t_2}) \left( \tilde{\gamma}_{s_1 m}
  \overset{l}{\circ} \left( \tilde{\gamma}_{m t_1} \overset{p}{\ast}
  \tilde{\gamma}_{s_2 t_2} \right) \right) \left( \tilde{\gamma}_{s_2 t_2}
  \overset{p}{\ast} \tilde{\gamma}_{m t_1} \right)\\
  &  & - \left\{\begin{array}{ll}
    \tilde{\gamma}_{s_1 m} \left( \tilde{\gamma}_{s_2 t_2} \overset{l}{\circ}
    \tilde{\gamma}_{m t_1} \right), & t_2 = m, \theta_m > \theta_t^{(2)}\\
    0, & \text{etc.}
  \end{array}\right.\\
  &  & + \left\{\begin{array}{ll}
    \left( \tilde{\gamma}_{s_1 m} \overset{l}{\circ} \tilde{\gamma}_{s_2 t_2}
    \right) \tilde{\gamma}_{m t_1}, & m = s_2, \theta_s^{(2)} > \theta_m\\
    0, & \text{etc.}
  \end{array}\right.\\
  &  & + \delta_{t_2 s_1} \begin{array}{c}
    \delta_{\theta_s^{(1)} < \theta_t^{(2)}} \tilde{\gamma}_{s_2 t_2}
    \overset{r}{\circ} \tilde{\gamma}_{s_1 m} \overset{l}{\circ}
    \tilde{\gamma}_{m t_1}\\
    \delta_{\theta_s^{(1)} > \theta_t^{(2)}} \tilde{\gamma}_{s_2 t_2}
    \overset{l}{\circ} \tilde{\gamma}_{s_1 m} \overset{l}{\circ}
    \tilde{\gamma}_{m t_1}
  \end{array} - \delta_{t_1 s_2} \begin{array}{c}
    \delta_{\theta_s^{(2)} < \theta_t^{(1)}} \tilde{\gamma}_{s_1 m}
    \overset{l}{\circ} \tilde{\gamma}_{m t_1} \overset{r}{\circ}
    \tilde{\gamma}_{s_2 t_2}\\
    \delta_{\theta_s^{(2)} > \theta_t^{(1)}} \tilde{\gamma}_{s_1 m}
    \overset{l}{\circ} \tilde{\gamma}_{m t_1} \overset{l}{\circ}
    \tilde{\gamma}_{s_2 t_2}
  \end{array} .
\end{eqnarray*}
Thus we have
\begin{eqnarray*}
  (\tilde{\gamma}_{s_1 m} \tilde{\gamma}_{m t_1}) \wedge \tilde{\gamma}_{s_2
  t_2} & = & \left( \sum_{p \in \gamma_{s_1 m} \cap \gamma_{s_2 t_2}}
  \tmop{sgn}_p (\gamma_{s_1 m}, \gamma_{s_2 t_2}) \left( \tilde{\gamma}_{s_1
  m} \overset{p}{\ast} \tilde{\gamma}_{s_2 t_2} \right) \left(
  \tilde{\gamma}_{s_2 t_2} \overset{p}{\ast} \tilde{\gamma}_{s_1 m} \right)
  \right) \tilde{\gamma}_{m t_1}\\
  &  & + \tilde{\gamma}_{s_1 m} \left( \sum_{p \in \gamma_{m t_1} \cap
  \gamma_{s_2 t_2}} \tmop{sgn}_p (\gamma_{m t_1}, \gamma_{s_2 t_2}) \left(
  \tilde{\gamma}_{m t_1} \overset{p}{\ast} \tilde{\gamma}_{s_2 t_2} \right)
  \left( \tilde{\gamma}_{s_2 t_2} \overset{p}{\ast} \tilde{\gamma}_{m t_1}
  \right) \right)\\
  &  & - \left( \delta_{m s_2} \begin{array}{c}
    \delta_{\theta_s^{(2)} < \theta_m} \tilde{\gamma}_{s_1 m}
    \overset{r}{\circ} \tilde{\gamma}_{s_2 t_2}\\
    \delta_{\theta_s^{(2)} > \theta_m} \tilde{\gamma}_{s_1 m}
    \overset{l}{\circ} \tilde{\gamma}_{s_2 t_2}
  \end{array} \right) \tilde{\gamma}_{m t_1} + \tilde{\gamma}_{s_1 m} \left(
  \delta_{t_2 m} \begin{array}{cc}
    \delta_{\theta_m < \theta_t^{(2)}} \tilde{\gamma}_{s_2 t_2}
    \overset{r}{\circ} \gamma_{12} & \\
    \delta_{\theta_m > \theta_t^{(2)}} \tilde{\gamma}_{s_2 t_2}
    \overset{l}{\circ} \gamma_{12} & 
  \end{array} \right)\\
  &  & + \left( \delta_{t_2 s_1} \begin{array}{c}
    \delta_{\theta_s^{(1)} < \theta_t^{(2)}} \tilde{\gamma}_{s_2 t_2}
    \overset{r}{\circ} \tilde{\gamma}_{s_1 m}\\
    \delta_{\theta_s^{(1)} > \theta_t^{(2)}} \tilde{\gamma}_{s_2 t_2}
    \overset{l}{\circ} \tilde{\gamma}_{s_1 m}
  \end{array} \right) \tilde{\gamma}_{m t_1} - \tilde{\gamma}_{s_1 m} \left(
  \delta_{t_1 s_2} \begin{array}{c}
    \delta_{\theta_s^{(2)} < \theta_t^{(1)}} \tilde{\gamma}_{m t_1}
    \overset{r}{\circ} \tilde{\gamma}_{s_2 t_2}\\
    \delta_{\theta_s^{(2)} > \theta_t^{(1)}} \tilde{\gamma}_{m t_1}
    \overset{l}{\circ} \tilde{\gamma}_{s_2 t_2}
  \end{array} \right)\\
  & = & (\tilde{\gamma}_{s_1 m} \wedge \tilde{\gamma}_{s_2 t_2})
  \tilde{\gamma}_{m t_1} + \tilde{\gamma}_{s_1 m} (\tilde{\gamma}_{m t_1}
  \wedge \tilde{\gamma}_{s_2 t_2}),
\end{eqnarray*}
and we finish the proof.
\end{prf}

\begin{pro}\label{Pro:PathtoUq}
    Suppose that the representative element $\tilde{\gamma}_{s_2t_2}$
    of path class with arguments is not self-intersecting, 
    particularly with $s_2\neq t_2$.
    Under the morphsim 
    $\mathfrak{S}(u_1,\ldots,u_\nu)
    \xrightarrow{\pi}
    U_q(\mathfrak{gl}_\nu)$ 
    induced by Theorem \ref{Thm:SWeight} and
    Theorem \ref{qStokes}, 
    we have
    \begin{align}\label{PartComRel}
        \pi(\tilde{\gamma}_{s_1t_1} \wedge
        \tilde{\gamma}_{s_2t_2}) =
        \frac{ 
        q^{-t( \tilde{\gamma}_{s_1t_1}, \tilde{\gamma}_{s_2t_2} )}
        \pi( \tilde{\gamma}_{s_1t_1} )
        \pi( \tilde{\gamma}_{s_2t_2} ) -
        q^{s( \tilde{\gamma}_{s_1t_1}, \tilde{\gamma}_{s_2t_2} )}
        \pi( \tilde{\gamma}_{s_2t_2} )
        \pi( \tilde{\gamma}_{s_1t_1} )}
        {q-q^{-1} },
    \end{align}
    where 
    $| \theta_s^{(2)} - \theta_s^{(1)} |,
     | \theta_t^{(2)} - \theta_s^{(1)} |,
     | \theta_t^{(2)} - \theta_t^{(1)} |,
     | \theta_s^{(2)} - \theta_t^{(1)} | < \pi$, and
    \begin{align}
        s( \tilde{\gamma}_{s_1t_1}, \tilde{\gamma}_{s_2t_2} ) & \assign
        \delta_{s_1s_2} \tmop{sgn}(\theta_s^{(2)} - \theta_s^{(1)}) -
        \delta_{s_1t_2} \tmop{sgn}(\theta_t^{(2)} - \theta_s^{(1)}),\\
        t( \tilde{\gamma}_{s_1t_1}, \tilde{\gamma}_{s_2t_2} ) & \assign
        \delta_{t_1t_2} \tmop{sgn}(\theta_t^{(2)} - \theta_t^{(1)}) -
        \delta_{t_1s_2} \tmop{sgn}(\theta_s^{(2)} - \theta_t^{(1)}).
    \end{align}
\end{pro}
\begin{prf}
    By analytically moving the endpoint of $\tilde{\gamma}_{s_2t_2}$,
    we can assume that it is a representative element of $\boldsymbol{\gamma}_{s_2t_2}^{(d)}$ without loss of generality.
    Subsequently, we only need to note that the right-hand side of equation \eqref{PartComRel} also satisfies the right Leibniz law.
    From Lemma \ref{Lem:Leb},
    it is suffices to verify the case where $\tilde{\gamma}_{s_1t_1}$ is also a representative element of
    $\boldsymbol{\gamma}_{s_1t_1}^{(d)}$ or
    $\boldsymbol{e}_{s}^{l}$, $\boldsymbol{e}_{s}^{r}$, 
    which reduces to equations \eqref{ComRelBegin}-\eqref{ComRelEnd}
    given by Theorem \ref{qStokes}.

    It can be checked that once the local behavior of 
    $\gamma_{s_1t_1}$ and $\gamma_{s_2t_2}$
    at their endpoints is fixed, 
    the wedge
    defined in \eqref{wedge} remains invariant under homotopy.
    Thus, equation \eqref{PartComRel} actually provides the commutation relation between a path class with arguments $\boldsymbol{\gamma}_{s_2t_2}$ with a representative element that does not self-intersect and any path class with arguments $\boldsymbol{\gamma}_{s_1t_1}$. 
\end{prf}

Based on this proposition, 
we can also derive the commutation relations between any two path classes with arguments. However, in this paper, we will not explore them in further detail.

\subsection{Braid Group Actions on Path Algebroids and Stokes Matrices}

\begin{figure}[h!]
    \centering
    \includegraphics[scale=1]{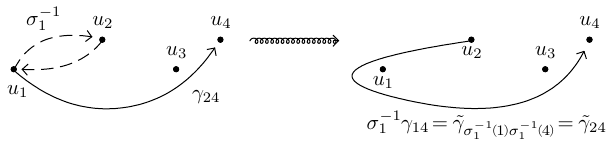}
    \caption{Example of the action of the element $\sigma_1$ 
    in the braid group $B_4$
    \label{Fig:test}}
\end{figure}

Let $B_\nu$ be the braid group with $\nu$ strands. For fixed $u_1,\ldots,u_\nu\in\mathbb{C}$, 
each $\sigma\in B_\nu$ induces an isomorphism between two path algebroids
\begin{align}\label{BraidAct}
\mathfrak{S}(u_1,\ldots,u_\nu,\infty)
\xrightarrow[\cong]{\sigma}
\mathfrak{S}(u_{\sigma 1},\ldots,u_{\sigma \nu},\infty).
\end{align}
See Figure \ref{Fig:test} for an illustration.


\begin{pro}\label{Def:sigmaS}
    Suppose that $(\mathcal{S},V_1,\ldots,V_\nu)$ is a representation of $\mathfrak{S}(u_1,\ldots,u_\nu,\infty)$,
    and $\sigma$ is an element of the braid group $B_\nu$.
    For every path class with arguments $\boldsymbol{\gamma}$, we denote
    \begin{align}\label{sigmaS}
    (\sigma \mathcal{S})(\boldsymbol{\gamma})
    & \assign 
    \mathcal{S}(\sigma^{-1}\boldsymbol{\gamma}).
    \end{align}
Then,
$\sigma(\mathcal{S},V_1,\ldots,V_\nu) \assign
(\sigma\mathcal{S},V_{\sigma^{-1}(1)},\ldots,V_{\sigma^{-1}(\nu)})$ is a representation of $\mathfrak{S}(u_1,\ldots,u_\nu,\infty)$.
\end{pro}


In the rest of this section, let us take any fixed direction $d$ such that
    \begin{equation*}
      \tmop{Im} (u_{1} \mathe^{\mathi d}) >
      \tmop{Im} (u_{2} \mathe^{\mathi d}) > \cdots > 
      \tmop{Im} (u_{\nu}  \mathe^{\mathi d}).
    \end{equation*}    
We define the generators \(\sigma_i\) of the braid group $B_\nu$ that interchange \(u_i\) and \(u_{i+1}\) by moving \(u_i\) anticlockwise around \(u_{i+1}\)
(during this process, we require that the projections of \(u_i\) and \(u_{i+1}\) in the \(-d\) direction do not pass through the projections of the other \(u_k\)'s).
On the one hand, the representation $(\mathcal{S},V_1,\ldots,V_\nu)$ can be encoded into a
pair of upper/lower triangular $(n_1,\ldots, n_\nu)$-block matrices $X^\pm=(X^\pm_{ij})_{\nu\times\nu}$ with entries
\begin{equation}\label{upperlower}
    X_{st}^+\assign
    \left\{\begin{array}{ll}
    \mathcal{S}(\boldsymbol{\gamma}_{st}^{(d)}) & ; s < t\\
    I & ; s = t\\
    0 & ; s > t
    \end{array}\right.,\quad
    X_{st}^-\assign
    \left\{\begin{array}{ll}
    0 & ; s < t\\
    I & ; s = t\\
    -\mathcal{S}(\boldsymbol{\gamma}_{st}^{(d)}) & ; s > t
    \end{array}\right..
    \end{equation}
Here recall that the paths $\boldsymbol{\gamma}_{st}^{(d)}$ are defined in Proposition \ref{Pro:gammajk}, and each entry $X_{st}^\pm$ has size $n_s\times n_t$.
On the other hand, the new representation
$(\sigma_i\mathcal{S},V_{\sigma_i^{-1}(1)},\ldots,V_{\sigma_i^{-1}(\nu)})$ gives rise to another pair of upper/lower triangular
$(n_{\sigma_i^{-1}(1)},\ldots, n_{\sigma_i^{-1}(\nu)})$-block matrices
$(\tilde{X}^+,\tilde{X}^-)$ with
\begin{align}\label{upperlowerS}
    \tilde{X}_{st}^+\assign
    \left\{\begin{array}{ll}
    (\sigma_i \mathcal{S})(\boldsymbol{\gamma}_{st}^{(d)}) & ; s < t\\
    I & ; s = t\\
    0 & ; s > t
    \end{array}\right.,\quad
    \tilde{X}_{st}^-\assign
    \left\{\begin{array}{ll}
    0 & ; s < t\\
    I & ; s = t\\
    -(\sigma_i \mathcal{S})(\boldsymbol{\gamma}_{st}^{(d)}) & ; s > t
    \end{array}\right..
\end{align}

In this way,
the braid group action of $\sigma$
on the representations induce 
an action on the pair of upper/lower triangular
$(n_{1},\ldots, n_{\nu})$-block matrices,
to produce a pair of upper/lower triangular
$(n_{\sigma^{-1}(1)},\ldots, n_{\sigma^{-1}(\nu)})$-block matrices,
by requiring
\[\sigma_i (X^+, X^-):=(\tilde{X}^+,\tilde{X}^-), \ \ \text{ for } \sigma_i\in B_\nu. \]
Here $\tilde{X}^\pm$ are given by \eqref{upperlowerS}.
Furthermore, by \eqref{upperlowerS} the entries of 
$\sigma_i(X^+,X^-)$ can be written as explicit polynomials of entries of $X^\pm$. It is determined by the corresponding topological relation.
For example, when $-d=\frac{\pi}{2}$,
Figure \ref{Fig:test} shows that the $(1,4)$-entry of the first component of 
$\sigma_1(X^+,X^-)$ is 
\[(\sigma_1\mathcal{S})(\boldsymbol{\gamma}_{14}^{(d)})=
\mathcal{S}(\sigma_1^{-1}\boldsymbol{\gamma}_{14}^{(d)})
=\mathcal{S}(\boldsymbol{\gamma}_{24}^{(d)})+
\mathcal{S}(\boldsymbol{\gamma}_{21}^{(d)})
\mathcal{S}(\boldsymbol{\gamma}_{14}^{(d)}).
\]
In general, we have the following result
\begin{cor}\label{expbraid}
    Suppose that $X^\pm=(X^\pm_{ij})_{\nu\times\nu}$ are 
    upper/lower triangular $(n_1,\ldots, n_\nu)$-block matrices,
    and all the diagonal blocks $X_{11},\ldots,X_{\nu\nu}$
    are identity matrices.
    Then the induced action of standard generator $\sigma_i (1\leqslant i \leqslant \nu-1)$ of braid group $B_\nu$
    on $(X^+,X^-)$ are
    pair of upper/lower triangular
    $(n_{\sigma_i^{-1}(1)},\ldots, 
    n_{\sigma_i^{-1}(\nu)})$-block matrices
    \begin{equation} \label{braid1}
    \sigma_i (X^+, X^-) =
    (K^-_i X^+ K_i^+, K^-_i X^- K_i^+),\quad
    1\leqslant i \leqslant \nu-1,
    \end{equation}
    where the block matrices $K^\pm_i$ are
    \begin{subequations}
    \begin{align}
    K^-_i & = 
    \tmop{diag} \left( I_{n_1}, \ldots, I_{n_{i - 1}},
    \left(\begin{array}{cc}
    - X^-_{i + 1, i} & I_{n_{i + 1}}\\
    I_{n_i} & O
    \end{array}\right), I_{n_{i + 2}}, \ldots, I_{n_{\nu}} \right),\\ \label{braid3}
    K_i^+ & = 
    \tmop{diag} \left( I_{n_1}, \ldots, I_{n_{i - 1}},
    \left(\begin{array}{cc}
    - X^+_{i, i + 1} & I_{n_i}\\
    I_{n_{i + 1}} & O
    \end{array}\right), I_{n_{i + 2}}, \ldots, I_{n_{\nu}} \right).
    \end{align}
    \end{subequations}
\end{cor}

If we consider the representations of the path algebroids from the Stokes phenomenon of the differential equations, then we obtain the braid group action on the space of Stokes matrices. 
When all multiplicities $n_1, ..., n_\nu$ equal to one, it recovers the well known braid group action on (non-blocking) Stokes matrices in the literature, see, e.g., \cite{Dubrovin, BoalchG}. Here we give a new interpretation via the natural braid group action on the path algebroids.

The braid group action on the space of Stokes matrices was first introduced to characterize the monodromy representation of the nonlinear isomonodromy deformation equations of the linear system \eqref{Conflu}. It is a differential equation
for the $(n_1,\ldots,\nu_\nu)$-block matrix valued function $A(u_1,...,u_\nu)$ with respect to the complex variables $u_1,...,u_\nu$:
\begin{eqnarray}\label{isoeq}
\frac{\partial}{\partial u_k}
A(u)=
[{\rm ad}^{-1}_{u}
{\rm ad}_{\mathcal{E}_{k}}A(u)
,A(u)],
\quad 
\text{for every $k=1,\ldots,\nu$},
\end{eqnarray} 
where
\begin{itemize}
    \item $\mathcal{E}_k=\tmop{diag} (0, \ldots ,
    \underset{k \text{-th}}{I_{n_k}} 
    , \ldots, 0)$
    as an $(n_1,\ldots,n_\nu)$-block diagonal matrix;

    \item For any $\nu\times \nu$ block matrix 
    $(X_{ij})_{\nu\times\nu}$, 
    denote $\ad_u^{-1} X$ as an
    $\nu\times \nu$ block matrix with the $(i,j)$-block entry
    \begin{align}\label{Nota:adu-1}
    (\ad_u^{-1} X)_{ij}\assign
    \left\{\begin{array}{ll}
    \frac{1}{u_i - u_j} X_{i j} & ; i \neq j\\
    0 & ; i = j
    \end{array}\right..
    \end{align}
\end{itemize}
The equation \eqref{isoeq} is deduced by the conditions that the Stokes matrices of the system \eqref{Conflu} with residue matrix $A(u)$ are locally constant while varying $u$. See \cite{JMMS} for more details. 
Following Miwa \cite{Miwa}, the solutions $A(u)$ with $u_1,...,u_\nu\in \mathbb{C}$ have the strong Painlev\'{e} property: they are multi-valued meromorphic functions of $u_1,...,u_\nu$ and the branching occurs when $u$ moves along a loop around the fat diagonal
\[\Delta=\{(u_1,...,u_\nu)\in \mathbb{C}^\nu~|~u_i = u_j, \text{for some } i\ne j \}.\]
Under the Riemann-Hilbert-Birkhoff map, the monodromy of the solutions $A(u)$ after the continuation along paths on $\mathbb{C}^\nu\setminus \Delta$ can be described geometrically as the braid group action on the corresponding Stokes matrices.
\begin{pro}\label{Prop:sigmaS}
   Suppose that $A(u_1,\ldots,u_\nu)$ satisfies 
   the isomonodromy equation \eqref{isoeq}
   corresponding to the system \eqref{Conflu}
   and $\sigma$ is an element of the braid group $B_\nu$.
   Denote
    \begin{subequations}
    \begin{align}
       \sigma u & \assign 
       \diag(u_1 I_{n_{\sigma^{-1}(1)}},\ldots,
       u_\nu I_{n_{\sigma^{-1}(\nu)}}),\\
       (\sigma A)(u_{\sigma(1)},\ldots,u_{\sigma(\nu)}) & \assign
       P_\sigma\cdot A(\sigma\curvearrowright(u_1,\ldots,u_\nu)) \cdot P_\sigma^{-1},
    \end{align}
    \end{subequations}
    where $A(\sigma\curvearrowright(u_1,\ldots,u_\nu))$ is the value of $A(u)$ after continuing from $(u_1,\ldots,u_\nu)$ to $(u_{\sigma(1)},\ldots,u_{\sigma(\nu)})$ along $\sigma$,
    and $P_\sigma=(\delta_{i\sigma(j)}I_{n_j})_{\nu\times\nu}$ is the permutation matrix.
    As a representation of $\mathfrak{S}(u_1,\ldots,u_\nu)$, we have
    \begin{align}\label{sigmauA}
        \sigma \mathcal{S}_{u,A} & =
        \mathcal{S}_{\sigma u,\sigma A},\\
        \label{BraidonS}
        \sigma(S_d^+(u,A),S_d^-(u,A)) & =
        (S_d^+(\sigma u, \sigma A),S_d^-(\sigma u, \sigma A)).
    \end{align}
\end{pro}

\begin{prf}
One first verifies that 
$\sigma \mathcal{S}_{u,A}, 
\mathcal{S}_{\sigma u,\sigma A}$ are representations of $\mathfrak{S}(u_1,\ldots,u_\nu)$,
\[\sigma \mathcal{S}_{u,A}(\boldsymbol{e}_s^l) =
\mathcal{S}_{\sigma u,\sigma A}(\boldsymbol{e}_s^l),\quad
s=1,\ldots,\nu.
\]
According to Remark \ref{Rmk:SW}, 
we only need to verify that 
for any distinct $s,t\in\{1,\ldots,\nu\}$,
we have 
$(\sigma \mathcal{S}_{u,A}) (\boldsymbol\gamma_{st}^{(d)}) =
\mathcal{S}_{\sigma u,\sigma A} (\boldsymbol\gamma_{st}^{(d)})$,
which can be guaranteed by the isomonodromy property.
Then applying \eqref{sigmauA} in Proposition \ref{Prop:sigmaS} and \eqref{sigmaS} in Definition \ref{Def:sigmaS}, one verifies \eqref{BraidonS}.
\end{prf}

\begin{rmk}
One can also consider the isomonodromy deformation of the quantum confluent hypergeometric equation \eqref{introqeq}, i.e., the equation
\begin{equation}
\frac{\mathd}{\mathd z} F(z) =
\left( u - \hbar \mathbf{E}(u) z^{-1}\right)\cdot F(z),
\end{equation}
where $\mathbf{E}(u)$ has the initial value $\mathbf{E}(u^0)=\mathbf{E}$ at a base point $u^0$, and is such that
the Stokes matrices are locally constant. Then the matrix function $E(u)$ satisfies the nonliear equation of the form \eqref{isoeq}. Applying Proposition \ref{Prop:sigmaS} to this case, while taking the base point $u^0$, leads to a braid group action on the Stokes matrices of the equation \eqref{introqeq}. Together with Theorem \ref{qStokes} or \ref{Squantumgroup}, one gets an action of $\sigma\in B_\nu$ on the generators of quantum group $U_q(\mathfrak{gl}_\nu)$. One checks by the formula in Corollary \ref{expbraid} that the braid group action coincides with the action, introduced by 
\cite{lusztig1990quantum}m
and independently by 
\cite{kirillov1990q} and 
\cite{soibel1990algebra}.
\end{rmk}

\section{Quantum Stokes Matrices and Yangian}\label{Sect:Q}



Having the preliminary given in the previous sections, we give proofs of Theorem \ref{MainThm}, Proposition \ref{infprodS} and Theorem \ref{oneside} in this section. 

\subsection{Proof of Theorem \ref{MainThm} and Proposition \ref{infprodS}}

\begin{prf}[Proofs of Theorem \ref{MainThm}]
Under assumption $\hbar\notin \mathbb{Q}$, the existence and uniqueness of the formal solution follow from Proposition \ref{FormSolirr}. But still to illustrate the process, we give a proof for our particular system. 

{\bf The existence and uniqueness of the formal solution.}
Suppose that a representation $\rho:U(\mathfrak{gl}_\nu)\rightarrow {\rm End}(V)$ is given.
Let us first assume the system \eqref{introqeq} associated to $V$, 
\begin{equation}\label{mn}
    \frac{\mathd}{\mathd z} {F}(z) =
    \left(u-\hbar\mathbf{E}^V z^{-1}\right)\cdot {F}(z),
\end{equation}
has a formal solution of the form given in \eqref{SolA}, i.e., 
\[
\hat{F}(z)= \left(I+\sum_{p=1}^\infty {H}^V_p z^{-p}\right)\cdot 
\e^{u z}
z^{-\hbar\delta\mathbf{E}^V},\quad
\delta\mathbf{E}^V=\diag(e^V_{11},\ldots,e^V_{\nu\nu}).\]
Here we denote $X^V:=\rho(X)$ for any $X\in U(\mathfrak{gl}_\nu)$, and each $H_p^V$ is an $\nu\times \nu$ matrix with entries valued in ${\rm End}(V)$.

Then glugging \eqref{SolA} into the equation \eqref{introqeq} and
comparing the coefficients of $z^{-m-1}$, we see that $H^V_m$ satisfies 
\begin{eqnarray}\label{simHm}
[u, H^V_{m+1}]=(-{m}+\hbar\mathbf{E}^V)\cdot H^V_{m}-H_{m}^V\cdot   \hbar\delta \mathbf{E}^V\in {\rm End}(V)\otimes {\rm End}(\mathbb{C}^\nu).\end{eqnarray}
Let $\{E_{kl}\}_{1\le k,l\le n}$ be the standard basis of ${\rm End}(\mathbb{C}^\nu)$. Recall that $\mathbf{E}=\sum_{k,l} e_{kl}\otimes E_{kl},$ and $u= \sum_i 1\otimes u_iE_{ii}$. Plugging $H^V_m=\sum_{k,l} H_{m, kl}\otimes E_{kl}$, with $H_{m,kl}\in {\rm End}(V)$, into the equation \eqref{simHm} gives rise to 
\begin{align} \nonumber
&\sum_{k,l}\frac{u_k-u_l}{\hbar} H_{m+1, kl}\otimes E_{kl}\\ \label{recuH}
=&-\sum_{k,l}\frac{m}{\hbar}H_{m,kl}\otimes E_{kl}+
\sum_{k,l,j} e^V_{kj} H_{m, jl} \otimes E_{kl}-\sum_{k,l}  H_{m, kl} e^V_{ll}\otimes E_{kl}.
\end{align}
That is for $k\ne l$ 
\begin{equation}\label{Hknel}
\frac{u_k-u_l}{\hbar} H_{m+1, kl}=-\frac{m}{\hbar}H_{m,kl}+\sum_{j=1}^n e^V_{kj} H_{m, jl}- H_{m, kl} e^V_{ll} \ \in {\rm End}(V),
\end{equation}
and for $k= l$ (replacing $m$ by $m+1$ in \eqref{recuH}), 
\begin{eqnarray}\label{Hkel}
0=\sum_{j\ne k} e^V_{kj} H_{m+1, jk}-\frac{m+1}{\hbar} H_{m+1, kk}+[e^V_{kk}, H_{m+1, kk}]\ \in {\rm End}(V).
\end{eqnarray}

Suppose $H^V_m$ is given, let us check that the above recursive relation has a unique solution $H^V_{m+1}$. First note that, since $u_k\ne u_l$ for $k\ne l$, the identity \eqref{Hknel} uniquely defines the "off-diangonal" part $H_{m+1, kl}$ ($k\ne l$) of $H_{m+1}$ from $H_m$. Furthermore, due to the assumption $\hbar\notin \mathbb{Q}$, we have $-\frac{m+1}{\hbar}{\rm Id}+{\rm ad}_{e^V_{kk}}$ is invertible on ${\rm End}(V)$ for any integer $m+1$. Thus, the condition \eqref{Hkel} uniquely defines the "diagonal" part $H_{m+1, kk}$ of $H^V_{m+1}$ from the off diagonal part.

{\bf The explicit form of $H_p$.} Note that the system \eqref{introqeq} associated to the representation $V$ becomes a special case of the System \eqref{Conflu} of rank $n=m\nu$, where $u=\diag(u_1 I_{m},\ldots, u_\nu I_{m})$, with distinct $u_1,\ldots,u_\nu$ and multiplicity $m,\ldots,m$, and
 the residue $A=-\hbar\mathbf{E}^V=-(\hbar e^V_{i j})_{\nu\times\nu}$ is divided into $(m,\ldots,m)$-blocks.

Therefore, we can apply Proposition \ref{FormSolirr} to the system \eqref{mn}. In particular, from \eqref{Cof:Finfhatkk}, \eqref{Cof:Finfkk} we have
\begin{subequations}
\begin{align}
    ({H}_{1}^V)_{\hat{k} k} & = 
    -\frac{1}{u_kI-u_{\hat{k}\hat{k}}}
    \hbar \mathbf{E}^V_{\hat{k}k},\\
    ({H}_{p+1}^V)_{\hat{k} k} & = 
    L_k(pI+\hbar (e_{kk}^V)^{\bm{r}}) \cdot
    ({H}_{p}^V)_{\hat{k} k},\\
    ({H}_p^V)_{kk} & =
    \frac{\hbar}{pI+\hbar(e_{kk}^V)^{\bm{r}}-\hbar e_{kk}^V} 
    \mathbf{E}_{k\hat{k}} \cdot
    ({H}_{p}^V)_{\hat{k} k}, 
\end{align}
\end{subequations}
where
\begin{align}\label{TypeALk}
    L_k(z) & =
    \frac{1}{u_k I - u_{\hat{k} \hat{k}}} 
    \left( z I - \hbar \mathbf{E}^V_{\hat{k} \hat{k}}- 
    \hbar^2  
    \frac{(\mathbf{E}_{\hat{k} k}\mathbf{E}_{k \hat{k}})^V}
    {(z-\hbar)I-\hbar e^V_{k k}} 
    \right).  
\end{align}
Let us inductively assume that for $p\geqslant 1$, we have 
\begin{align}
\label{HVhatkkcomrel}
    (e_{k k}^V)^{\bm{r}} \cdot ({H}^V_{p})_{\hat{k}k} =
	({H}^V_{p})_{\hat{k}k} e_{k k}^V =
	\tmop{diag} (\underset{\text{$n - 1$ terms}}{\underbrace{e_{k k}^V + I,
    \ldots, e_{k k}^V + I}}) \cdot
    ({H}^V_{p})_{\hat{k} k},
\end{align}
thus we have
\begin{subequations}
\begin{align}\label{Eq:YangCoef}
    ({H}^V_{p+1})_{\hat{k}k} & =
    L_k((p+\hbar)I+\hbar e_{kk}^V) \cdot
    ({H}_{p}^V)_{\hat{k} k} =
    T_k(p)^V \cdot
    ({H}_{p}^V)_{\hat{k} k},\\
    \label{Eq:YangCoefb}
    ({H}^V_{p})_{k k} & =
    \frac{\hbar \mathbf{E}^V_{k \hat{k}}}{p} 
    ({H}^V_{p})_{\hat{k} k}.
\end{align}
\end{subequations}
It can be seen that the base case $p=1$ holds for \eqref{HVhatkkcomrel}, and $T_k(p)$ commutes with $\tmop{diag} (e_{k k}^V + I,
 \ldots, e_{k k}^V + I)$.
Thus, \eqref{HVhatkkcomrel} also holds for $p+1$, and induction ensures that it holds for any $p\geqslant 1$.
Therefore, \eqref{Eq:YangCoef}, \eqref{Eq:YangCoefb} also holds for any $p\geqslant 1$.

The above argument works for any $V$.
Since the representations of $U(\mathfrak{gl}_\nu{})$ separates the elements in $U(\mathfrak{gl}_\nu{})$, we have \eqref{Form:QHp}.

{\bf Yangian relations.}
In the end, the relation \eqref{Intro:RTTYangian} can be directly verified. There is another way to check it: our $T_k(\lambda)$, arising from the formal solution of the quantum confluent hypergeometric equation \eqref{introqeq}, defines a $(u_1,...,u_\nu)$-family deformation $O(u)_k$ (as in \eqref{Oconstruct1}) of the map for each $k=1,...,\nu$
\begin{equation*}
 O_k: Y_{\hbar}(\mathfrak{gl}_{\nu-1})\rightarrow U(\mathfrak{gl}_\nu)~;~ T(\lambda)_{ij}\mapsto  
    (\lambda+\hbar(e_{kk}+1))\delta_{ij}
    - \hbar e_{ij}
    - \hbar^2\frac{e_{ik}e_{kj}}{\lambda}
    \in U(\mathfrak{gl}_\nu),\quad
    i,j\in \{ 1,\ldots,\nu\}\setminus \{k\}.
 \end{equation*}
The map $O_k$ is known as the Olshanski centralizer construction and is known to be an algebra homomorphism.
See e.g., in \cite[Theorem 1.12.1]{molev2007yangians} for more details on the algebra homomorphism. Then one checks that the deformations $O(u)_k$ are also algebra homomorphisms.
\end{prf}

\begin{prf}[Proof of Proposition \ref{infprodS}]

Suppose that a representation $\rho:U(\mathfrak{gl}_\nu)\rightarrow {\rm End}(V)$ is given. Let us consider the system \eqref{introqeq} associated to $V$. Such a system becomes a special case of the System \eqref{Conflu}, i.e., an $mn\times mn$ system (assume $m={\rm dim}(V)$), where $u=\diag(u_1 I_{m},\ldots, u_\nu I_{m})$, with distinct $u_1,\ldots,u_\nu$ and multiplicity $m,\ldots,m$, and
 the residue $A=-\hbar\mathbf{E}^V=-(\hbar e^V_{i j})_{\nu\times\nu}$ is divided into $(m,\ldots,m)$-blocks.

Therefore, we can apply 
the results of Theorem \ref{MainThm} and
Theorem \ref{Thm:Difeq} to the system \eqref{introqeq} associated to $V$. It gives rise to the desired expressions in Prosition \ref{infprodS}.
\end{prf}


\subsection{Proof of Theorem \ref{oneside}}

\begin{prf}[Proof of Theorem \ref{oneside}]
The matrix $\mathscr{T}$ in Theorem \ref{oneside}
represents the infinite product of matrices $T_k(m)$,
i.e.
\begin{align}\label{FormInfMatProd}
    \left( \overleftarrow{\prod_{m=1}^{p}}
	T_k(m)
    \right)
    \left(
    -\frac{1}{u_kI-u_{ \hat{k} \hat{k} }}
    \hbar \mathbf{E}_{\hat{k}k}
    \right) \overset{p\to\infty}{\sim}
    \mathscr{T}_{\hat{k} k},\quad
    \left(-\hbar \mathbf{E}_{k\hat{k}} \right)
    \left( \overleftarrow{\prod_{m=-p}^{-1}}
	T_k(m)
	\right) \overset{p\to\infty}{\sim}
    \mathscr{T}_{ k \hat{k}}.
\end{align}
Thus by Proposition \ref{infprodS} we have
\begin{align}
    \mathscr{T}_{st} \assign
    (u_t-u_s)^{\hbar e_{ss}}
    \frac{(S_{[\tau]})_{st}}{2\pi\mathi}
    (u_t-u_s)^{-\hbar e_{tt}},
    \quad s\neq t.
\end{align}
Thus the commutation relations of the entries of $\mathscr{T}$
can be reduced to the commutation relations of
the weights $\mathcal{S}$ of segments under
quantum system \eqref{introqeq}, 
as guaranteed by Proposition \ref{Pro:PathtoUq}.

The segments $[u_s,u_t]$ with argument,
as the path class with arguments, 
under the morphsim
$\mathfrak{S}(u_1,\ldots,u_\nu)
 \xrightarrow{\pi}
 U_q(\mathfrak{gl}_\nu{})$
induced by Theorem \ref{qStokes},
are mapped to
\begin{align*}
    [u_s,u_t] \mapsto
    \tilde{\mathscr{T}}_{st } \assign
    (u_t-u_s)^{-\hbar e_{ss}} \cdot
    2\pi\mathi \mathscr{T}_{st} \cdot
    (u_t-u_s)^{ \hbar e_{tt}}.
\end{align*}
According to Proposition \ref{Pro:PathtoUq},
for distinct indices $s_1,t_1,s_2,t_2$ and
$\arg(u_{t_2}-u_{s_2}) -
 \arg(u_{t_1}-u_{s_1}) \in (0,\pi)$, we have
\begin{subequations}
\begin{alignat}{2}
    \tilde{\mathscr{T}}_{s_1 t_1} 
    \tilde{\mathscr{T}}_{s_2 t_2} -
    \tilde{\mathscr{T}}_{s_2 t_2} 
    \tilde{\mathscr{T}}_{s_1 t_1} & = 
    0,& \quad
    [u_{s_1},u_{t_1}]\cap [u_{s_2},u_{t_2}] & = \varnothing,\\
    \frac{
    \tilde{\mathscr{T}}_{s_1 t_1} 
    \tilde{\mathscr{T}}_{s_2 t_2} -
    \tilde{\mathscr{T}}_{s_2 t_2} 
    \tilde{\mathscr{T}}_{s_1 t_1}
    }{q-q^{-1}} & = 
    \tilde{\mathscr{T}}_{s_1 t_2}
    \tilde{\mathscr{T}}_{s_2 t_1},& \quad
    [u_{s_1},u_{t_1}]\cap [u_{s_2},u_{t_2}] & \neq \varnothing.
\end{alignat}
\end{subequations}
For distinct indices $s,m,t$, we have
\begin{subequations}
\begin{align}
    \frac{q 
    \tilde{\mathscr{T}}_{sm}
    \tilde{\mathscr{T}}_{mt} -
    \tilde{\mathscr{T}}_{mt}
    \tilde{\mathscr{T}}_{sm}
    }{ q - q^{-1} } & =
    - \tilde{\mathscr{T}}_{st},\quad
    \tmop{arg}(u_t-u_m) -
    \tmop{arg}(u_m-u_s) \in (0,\pi),\\
    \frac{q^{-1}
    \tilde{\mathscr{T}}_{sm}
    \tilde{\mathscr{T}}_{mt} -
    \tilde{\mathscr{T}}_{mt}
    \tilde{\mathscr{T}}_{sm}
    }{ q - q^{-1} } & =
    - \tilde{\mathscr{T}}_{st},\quad
    \tmop{arg}(u_t-u_m) -
    \tmop{arg}(u_m-u_s) \in (-\pi,0),\\
    \tilde{\mathscr{T}}_{ms}
    \tilde{\mathscr{T}}_{mt} -
    q \tilde{\mathscr{T}}_{mt}
    \tilde{\mathscr{T}}_{ms} & =
    0,\quad
    \tmop{arg}(u_t-u_m) -
    \tmop{arg}(u_s-u_m) \in (0,\pi),\\
    q^{-1} \tilde{\mathscr{T}}_{sm}
    \tilde{\mathscr{T}}_{tm} -
    \tilde{\mathscr{T}}_{tm}
    \tilde{\mathscr{T}}_{sm} & =
    0,\quad
    \tmop{arg}(u_m-u_t) -
    \tmop{arg}(u_m-u_s) \in (0,\pi),\\
    \frac{
    q \tilde{\mathscr{T}}_{st}
    \tilde{\mathscr{T}}_{ts} -
    q^{-1} \tilde{\mathscr{T}}_{ts}
    \tilde{\mathscr{T}}_{st}
    }{q-q^{-1}} & =
    q^{2(e_{ss}-e_{tt})} - 1,\quad
    \tmop{arg}(u_s-u_t) -
    \tmop{arg}(u_t-u_s) = \pi.
\end{align}
\end{subequations}
Applying
$q^{e_{kk}} \mathscr{T}_{st} q^{-e_{kk}} =
q^{\delta_{sk}-\delta_{kt}} \mathscr{T}_{st}$,
we finish the proof.
\end{prf}

\section*{Acknowledgements}
\noindent The authors are supported by the National Natural Science Foundation of China (No. 12171006) and by the National Key Research and Development Program of China (No. 2021YFA1002000).








\bibliography{paper.bib}
\bibliographystyle{plain}

\end{document}

%% file: preamble/usepackage.tex
\usepackage{geometry}
\usepackage{hyperref}
\usepackage{theoremref}

\usepackage{times}
\usepackage{amsthm}
\usepackage{scalerel}
\usepackage{amsfonts}
\usepackage{amssymb}
\usepackage{mathdots}
\usepackage{color}
\usepackage{mathrsfs}
\usepackage{multirow}

\usepackage{amsmath}
\usepackage{bm}
\usepackage{stmaryrd}
\SetSymbolFont{stmry}{bold}{U}{stmry}{m}{n}

\usepackage{tikz-cd}
\usepackage{tikz}
\usetikzlibrary{matrix,arrows,decorations.pathmorphing}
\usetikzlibrary{positioning}
\usepackage[all]{xy}
\usepackage{graphicx}
\usepackage{subfigure}

\usepackage[toc,page]{appendix}

\usepackage[usestackEOL]{stackengine}
\usepackage{dutchcal}

%% file: preamble/newcommand.tex
\newcommand{\ad}{\tmop{ad}}

\newcommand{\diag}{\tmop{diag}}

\newcommand{\e}{\mathrm{e}}

\newtheorem{thm}{Theorem}[section]
\newtheorem{lem}[thm]{Lemma}
\newtheorem{cor}[thm]{Corollary}
\newtheorem{pro}[thm]{Proposition}

\theoremstyle{definition}

\newtheorem{rmk}[thm]{Remark}
\newtheorem{defi}[thm]{Definition}
\newtheorem{nota}[thm]{Notation}
\numberwithin{equation}{section}

%% file: preamble/texmacs.tex
\usepackage[normalem]{ulem}

\newcommand{\assign}{:=}
\newcommand{\mathd}{\mathrm{d}}
\newcommand{\mathe}{\mathrm{e}}
\newcommand{\mathi}{\mathrm{i}}

\newcommand{\tmop}[1]{\ensuremath{\operatorname{#1}}}

\newcommand{\udots}{{\mathinner{
\mskip1mu\raise1pt\vbox{\kern7pt\hbox{.}}
\mskip2mu\raise4pt\hbox{.}
\mskip2mu\raise7pt\hbox{.}\mskip1mu}}}

\newenvironment{prf}
    {\proof}
    {\hspace*{\fill}$\Box$\medskip}
    